\documentclass[12pt]{amsart}
\usepackage{amsmath}
\numberwithin{equation}{section}
\usepackage{amsthm}
\usepackage{amssymb}
\usepackage{amsfonts}
\usepackage{pinlabel}
\usepackage{graphicx}
\usepackage{float}
\usepackage{relsize}
\newcommand{\nc}[2]{ \newcommand{#1}{#2} }

\nc{\avint}{ {- \hspace{-3.5mm} \int} } 

\nc{\R}{\mathrm{I \! R}} 
\nc{\N}{\mathrm{ I \! N}} 

\newcommand{\pclosure}[1]{ \stackrel{\rule{.5 in}{.01 in}}{#1} }

\DeclareRobustCommand{\rchi}{{\mathpalette\irchi\relax}}
\newcommand{\irchi}[2]{\raisebox{\depth}{$#1\chi$}}
\newcommand{\bigrchi}{\mathlarger{\mathlarger{\rchi}}}

\newcommand{\refeqn}[1]{ (\!\!~\ref{eq:#1}) } 
\newcommand{\refthm}[1]{ \!\!~\ref{#1} } 

\nc{\Holder}{H\"{o}lder\ }

\nc{\ith}{ \ensuremath{\text{i}^{\text{th}}} }
\nc{\jth}{ \ensuremath{\text{j}^{\text{th}}} }
\nc{\kth}{ \ensuremath{\text{k}^{\text{th}}} }
\nc{\dst}{ \ensuremath{\text{1}^{\text{st}}_{\delta}} }
\nc{\dnd}{ \ensuremath{\text{2}^{\text{nd}}_{\delta}} }
\nc{\ost}{ \ensuremath{\text{1}^{\text{st}}} }
\nc{\tnd}{ \ensuremath{\text{2}^{\text{nd}}} }
\nc{\curl}{ \nabla \times }
\nc{\Div}{ \nabla \cdot }
\nc{\DC}{K}

\nc{\Ppl}{ \mathcal{M}^{+} }  \nc{\Pmn}{ \mathcal{M}^{-} }

\nc{\smiley}{ $\stackrel{\because}{\smile} \;$ }

\renewcommand{\eqref}[1]{(\ref{#1})}
\newcommand{\eat}[1]{}

\newcommand{\bsmall}{\begin{array}[c]{c}}
\newcommand{\esmall}{\end{array}}

\renewcommand{\div}{\mathrm{div}}

\theoremstyle{plain}

\theoremstyle{definition}

\theoremstyle{plain}
\newtheorem{theorem}{Theorem}[section]
\newtheorem{mainthm}{Theorem}
\newtheorem{proposition}[theorem]{Proposition}
\newtheorem{lemma}[theorem]{Lemma}
\newtheorem{corollary}[theorem]{Corollary}

\theoremstyle{definition}
\newtheorem{definition}[theorem]{Definition}
\newtheorem{remark}[theorem]{Remark}

\def\qed{\hfill\rule{1ex}{1ex}\\}
\newenvironment{pf}{\noindent {\bf Proof.}}{\qed}

\numberwithin{equation}{section}

\def\>{ > _{\sigma}}

\title[Uniqueness and Counterexamples]{Strong and Weak Mean Value Sets
for Elliptic Operators: Uniqueness and Counterexamples}
\author[Armstrong]{Niles Armstrong}
\address{Milwaukee School of Engineering, Mathematics Department,
1025 N. Broadway,
Milwaukee, WI 53202}
 \email{armstrongn@msoe.edu}
\author[Blank]{Ivan Blank}
\address{Kansas State University, Mathematics Department,
138 Cardwell Hall,
Manhattan, KS 66506}
 \email{blanki@math.ksu.edu}

\begin{document}
\baselineskip 18pt

\begin{abstract}
Classical theorems of Epstein--Schiffer and Kuran
characterize balls by the mean value identity for
integrable harmonic functions. In the Fermi Lectures,
Caffarelli described an obstacle problem construction
of analogous averaging sets for symmetric uniformly elliptic
divergence-form operators. For such an operator with bounded
measurable coefficients in dimension $n \ge 2,$ fix a center
$x_0$ and a volume $V > 0,$ and let $E$ be the corresponding
canonical set. We prove that a bounded open set $D$ of volume
$V$ containing $x_0$ satisfies the mean value inequality at
$x_0$ for every integrable subsolution on $D$ if and only if
$E \subset D \subset \operatorname{int}\overline{E}.$
In particular, every such strong mean value set agrees with
$E$ almost everywhere. The proof uses weak modified Schwarz
potentials and an energy comparison.

Additionally, we construct smooth symmetric uniformly elliptic
operators for which the unit ball satisfies the mean value
identity for all integrable solutions but fails the subsolution
inequality, showing that there can be two distinct notions
for what one can consider to be a ``mean value set.''
The canonical mean value set with the same center and volume
contains the entire unit sphere and has a hole of positive
measure inside the unit ball.

\


\

\end{abstract}
\maketitle

\setcounter{tocdepth}{2}
\begingroup
\small
\setlength{\parskip}{0pt}
\setlength{\baselineskip}{12pt}
\tableofcontents
\endgroup

\section{Introduction}   \label{Intro}

Mean value properties raise an inverse geometric question:
which bounded open sets recover the value of a solution at
a prescribed point by an unweighted volume average?
For the Laplacian, the classical theorems of
Epstein--Schiffer and Kuran show that the mean value
identity for all integrable harmonic functions
characterizes balls \cite{EpSchiff, Ku}.
The corresponding mean value inequality for subharmonic
functions characterizes the same balls. Dalmasso showed
that it suffices to test this inequality on integrable
$C^2$ subharmonic functions that are not harmonic
\cite[Theorem~3]{Dal08}.
We study these two characterizations for symmetric
uniformly elliptic divergence-form operators and show
that they can lead to different classes of averaging sets.

Mean value formulas for uniformly elliptic divergence-form
operators with bounded measurable coefficients go back at
least to Littman, Stampacchia, and Weinberger in 1963
\cite[Section 8]{LSW}.
Their formula uses weighted averages of the solution over
regions defined by level sets of the Green's function,
with weights involving a quadratic form in its gradient.
For more recent weighted surface and volume mean value
formulas for uniformly elliptic divergence-form equations
with H\"older continuous coefficients, see
Pallara and Polidoro \cite{PP24}.
For a class of hypoelliptic divergence-form operators,
Abbondanza and Bonfiglioli show that a weighted mean value
identity characterizes the superlevel sets of the fundamental
solution with pole at the averaging point \cite{AbbBon13}.
In his 1998 Fermi lectures on the obstacle problem, Caffarelli
described, for $L := \div (A\nabla \cdot)$ defined on $\R^n$
(where $A$ is bounded, measurable, symmetric, and uniformly
elliptic), and for each point $x_0,$ a family of obstacle problems
whose unique solutions have nested noncontact sets $D_r(x_0)$
comparable to $B_r(x_0)$ \cite{C2}, and which give an
unweighted mean value inequality for the operator $L.$
In other words, for every subsolution $v$ of $Lv = 0,$ the unweighted
average of $v$ over $D_r(x_0)$ is nondecreasing in $r$.
We give the precise construction in Section\refthm{MVTOP}\!\!,
and we will refer to these sets as the ``canonical mean value
sets with center $x_0$ and volume $r^n.$'' (We will call this
the ``canonical set'' for short.)

Caffarelli and Roquejoffre \cite{CR} stated this mean value
theorem and referred to the book of Kinderlehrer and
Stampacchia \cite{KS}. As explained in \cite[p.~130]{BH1},
Blank, Hao, and Kinderlehrer were unable to locate in that
book a theorem directly establishing the existence of the
solution to the singular obstacle problem needed in the
proof. Blank and Hao subsequently supplied a detailed proof
of the mean value theorem in \cite{BH1}, including the
existence and nondegeneracy results needed for the
obstacle-problem construction. We have since studied
the geometry of the canonical sets $D_r(x_0)$
in \cite{Ar, ArmB, AryB, BBL, BH2}.
Throughout the paper, solutions and subsolutions used as mean
value test functions are understood in the usual weak sense
and belong to $W^{1,2}_{\mathrm{loc}}(U)$ on the open set $U$ where
the equation or inequality is imposed. When averaging such
a function over a mean value set $D\subset U$, we additionally
require that it belong to $L^1(D)$.
The equations for the singular Green functions and Schwarz
potentials are interpreted distributionally in the function
spaces specified below.

We start with a definition:
\begin{definition}[Weak mean value sets]    \label{MVsets}
A bounded open set $D$ is called a weak mean value set
for $L$ at $x_0 \in D$ if the mean value equality
\begin{equation}
     u(x_0) = \int_{D} \hspace{-.21in} \rule[.75mm]{.3cm}{.35mm} \ \ \ \
            u(x) \; dx 
\label{eq:MVprop}
\end{equation}
holds for every $u \in L^1(D)\cap W^{1,2}_{\mathrm{loc}}(D)$
satisfying $Lu = 0$ weakly in $D.$
Throughout the paper,
a slash through an integral denotes division by the
measure of its domain of integration.
We use the continuous representative of each weak solution.
For a weak subsolution $v,$ we use the upper semicontinuous
representative obtained by taking the negative of the lower
semicontinuous representative of $-v$ given by
\cite[Corollary~6.4]{BH1}.
\end{definition}

\begin{remark}[The domain of the test functions]
\label{Ambient}
The test functions in Definition~\refthm{MVsets}need only
be defined on $D.$ In particular, every weak mean value
set is connected: otherwise, the function equal to zero
on the component containing $x_0$ and equal to one on
the other components contradicts the mean value equality.

Checking the equality only for functions that extend
harmonically to a larger open set is not sufficient
in general. For example, for the Laplacian, the union
of a ball and a disjoint open concentric shell satisfies
the equality for every integrable harmonic function on
a concentric ball containing both components. Nevertheless,
the disconnected union is not a weak mean value set.
\end{remark}
\noindent
Requiring the corresponding inequality for subsolutions
gives the stronger notion below.
\begin{definition}[Strong mean value sets]    \label{SMVsets}
A bounded open set $D$ is called a strong mean value set
for $L$ at $x_0 \in D$ if the mean value inequality
\begin{equation}
     u(x_0) \le \int_{D} \hspace{-.21in} \rule[.75mm]{.3cm}{.35mm} \ \ \ \
            u(x) \; dx 
\label{eq:MVineq}
\end{equation}
holds for every $u \in L^1(D)\cap W^{1,2}_{\mathrm{loc}}(D)$
satisfying $Lu \ge 0$ weakly in $D.$
\end{definition}
\noindent
Every strong mean value set is a weak mean value set, as follows by
applying the subsolution inequality to both $u$ and
$-u$ when $Lu = 0.$

Epstein and Schiffer \cite{EpSchiff} established the following
classical characterization, stated here in our terminology.
Their hypotheses apply because our weak mean value sets are
bounded and connected; in particular, they have finite measure
and their complements have nonempty interior.

\begin{theorem}[Epstein--Schiffer, 1965] \label{KurThm}
If $D$ is a weak mean value set for the point $x_0$ and the
Laplacian, then $D$ is an open ball centered at $x_0.$
\end{theorem}
\noindent
See also Epstein's earlier result for simply connected
planar domains \cite{Ep} and the subsequent work of
Goldstein and Ow \cite{GO}. Kuran \cite{Ku} later
extended the characterization to arbitrary domains
of finite measure.
Cupini, Lanconelli, Fusco, and Zhong proved a sharp
quantitative stability version of this characterization
\cite[Theorem~1.1]{CFLZ}.
Their estimate controls the normalized symmetric
difference from the ball with the same center and
volume by a normalized defect in harmonic averaging.
Recent extensions of this inverse mean-value and
stability viewpoint include the sharp surface stability
results of Cupini and Lanconelli \cite{CL25}, the
polyharmonic rigidity and stability theorem of
Abatangelo \cite{Abat26}, and the converse mean-value
theorem of Morav\'ik and Youssfi for the invariant
Laplace--Beltrami operator \cite{MY26}.
We ask whether the classical characterization extends
to uniformly elliptic divergence-form operators when
balls are replaced by the canonical sets arising from
Caffarelli's obstacle problem.
Our results distinguish the roles of the mean value
identity and the subsolution inequality.

We first classify all open strong mean value sets with a given
center and volume. Two quantities enter the characterization.
For a weak mean value set $D$ at $x_0,$ choose a sufficiently
large ball $\Omega$ with $D \Subset \Omega,$ and let $G$ be the
positive Dirichlet Green function for $L$ in $\Omega.$
The weak modified Schwarz $L$-potential of $D$ is given, for
$z \in \Omega \setminus \{x_0\},$ by
$$W_D(z) := G(z,x_0) - \frac{1}{|D|}\int_D G(z,x)\,dx.$$
For $z \in \Omega \setminus D,$ the weak mean value identity
applied to $G(z,\cdot)$ gives $W_D(z) = 0$
(Theorem~\refthm{MVSYWMSP}\!\!).
The maximal mean-value gap $\delta(D)$ is the largest value
of
$$u(x_0) - |D|^{-1}\int_D u(x) \, dx$$
among solutions
$u \in W^{1,2}_0(\Omega)$ of $Lu = f$ with $0 \le f \le 1$
(Proposition~\ref{potdefect}).

\begin{mainthm}[Classification of strong mean value sets]
\label{WisSiffPisP}\label{openreps}
Let $n \ge 2,$ and let $L := \div(A\nabla\cdot)$ on $\R^n,$
where $A$ is bounded, measurable, symmetric, and uniformly
elliptic. Fix a center $x_0$ and a parameter $r > 0,$ and
set $E := D_r(x_0).$ Let $D$ be a bounded open set containing
$x_0,$ with $|D| = r^n.$ The following conditions are equivalent.
\begin{enumerate}
\item $D$ is a strong mean value set for $L$ at $x_0.$
\item $D$ is a weak mean value set for $L$ at $x_0,$ and
      its weak modified Schwarz $L$-potential is nonnegative.
\item $E \subset D \subset \operatorname{int}\overline{E}.$
\item $D$ is a weak mean value set for $L$ at $x_0,$ and
      its maximal mean-value gap $\delta(D)$ is zero.
\end{enumerate}
When these conditions hold, $|D \triangle E| = 0$ and
$W_D = W_E$ pointwise away from $x_0$ on any common bounded
auxiliary domain. Consequently, $E$ and
$\operatorname{int}\overline{E}$ are respectively the smallest
and largest open strong mean value sets with center $x_0$ and
volume $r^n.$ There is exactly one such open set if and only if
$E = \operatorname{int}\overline{E}.$
\end{mainthm}

For the Laplacian, $E$ is a ball and hence equals
$\operatorname{int}\overline{E}.$ Thus Theorem~\ref{openreps}
recovers the characterization of balls by the strong mean
value property. The distinction between open representatives
matters here because deleting even a single point can change
the admissible test functions and destroy the weak mean value
property (Proposition~\ref{pointdefect} and Remark~\ref{deletepoint}).

We next ask what the mean value identity alone determines.

\begin{mainthm}[Weak does not imply strong]
\label{WDNIS}
For every $n \ge 2,$ there exists a smooth, symmetric,
uniformly elliptic matrix field $A$ on $\R^n,$ with $A = I$
outside a compact subset of $D := B_1(0),$ such that, for
$L := \div(A\nabla\cdot)$ and center $0,$
$$D \text{ is a weak mean value set, but not a strong mean value set.}$$
When $n = 2,$ the matrix $A$ may additionally be chosen
with $\det A = 1.$
\end{mainthm}

The ball and the canonical set with the same center and volume
are both weak mean value sets and differ on a set of positive
measure. This disproves the uniqueness assertion for weak mean
value sets in our earlier preprint \cite{ABthepain}. The failure
occurs for a single positive point-mass quadrature datum, with
unweighted Lebesgue averages and the full class of integrable
solutions defined on the averaging set itself.
Let $E$ denote the canonical set with the same center and volume
as $D.$ The argument between (4.4) and (4.5) in \cite{ABthepain}
inferred from $|D \setminus E| > 0$ that some point of $\partial D$
lies outside $\overline{E}.$ Our third theorem shows that the
entire boundary of $D$ can lie inside $E.$

\begin{mainthm}[Strong mean value sets with a hole]
\label{ESMVWH}
Let $n \ge 2,$ and let $L := \div(A\nabla\cdot)$ on $\R^n,$
where $A$ is bounded, measurable, symmetric, and uniformly
elliptic. Let $D$ be a bounded weak mean value set for $L$
at $x_0$ which is not strong. Assume that $\partial D$ is $C^{1,1}$ and
that $A$ is smooth in an open neighborhood of $\partial D.$
Set $\mathcal{V} := |D|$ and $E := D_{\mathcal{V}^{1/n}}(x_0).$
Then $\partial D \subset E.$ Moreover, $E^c$ has a bounded
connected component of positive measure contained in $D$ and
an unbounded connected component of infinite measure.

In particular, these conclusions hold in every dimension
$n \ge 2$ for $D = B_1(0)$ and the globally smooth operator
of Theorem~\refthm{WDNIS}\!\!. The canonical set with the same center
and volume contains the entire unit sphere and misses an open
ball inside $B_1(0).$
\end{mainthm}

The almost-everywhere equality in Theorem~\refthm{openreps}is
only the measure-theoretic
part of our classification. Integral and variational
characterizations of sets naturally tend to determine
characteristic functions only up to null sets. This is the
usual level of uniqueness for subharmonic quadrature domains;
recent work distinguishes ``strong quadrature domains'' by
strict positivity of the associated potential, which restores
exact uniqueness \cite[Introduction]{KSS26}. A closer analogue
occurs for subharmonic balls: Shahgholian and Sj\"odin identify
smallest and largest open representatives between which every
subharmonic ball lies \cite[Remark 3.3]{ShSj}. In a different
setting, Bou-Rabee and Gwynne obtain exact uniqueness of
Liouville quantum gravity harmonic balls after requiring that
a ball equal the interior of its closure
\cite[Theorem 1.1]{BRG25}.

The broader potential-theoretic setting for our results is
the classical theory of quadrature domains. A mean value
identity is a quadrature identity with datum
$|D|\delta_{x_0}$; see Section~\ref{ACF}.
For the Laplacian, the construction of subharmonic
quadrature domains through partial balayage and obstacle
problems is classical; see
\cite[Sections 5--7]{GusShap}.
For point-mass data, Green-potential characterizations
and comparison results for harmonic and subharmonic balls
are given in
\cite[Theorem 2.1 and Propositions 3.1--3.2]{ShSj}.
Shahgholian and Sj\"odin ask whether a nontrivial harmonic
ball can fail to be a subharmonic ball and whether a
subharmonic ball can have a hole even when the containing
domain has none \cite[Questions~3.2 and~3.3]{ShSj}.
Their test functions are defined using the Green function
of the Laplacian in the containing domain.
Theorems~\refthm{WDNIS}and~\refthm{ESMVWH}exhibit both phenomena
for smooth variable coefficients on $\R^n,$ using all
integrable solutions and subsolutions on the averaging set
itself.

The connection between subharmonic quadrature inequalities
and partial balayage is also developed in the Riemannian
setting by Gustafsson and Roos
\cite[Sections 8 and 10]{GRooooos}.
Loutsenko used planar Hele--Shaw flow in an inhomogeneous
medium to derive quadrature identities for divergence-form
elliptic equations with a variable scalar coefficient
\cite{Loutsenko06}. For a single point source in an initially
empty cell with unit area weight, his identity says that the
area average of a solution over the evolving domain equals
its value at the source. His explicit disk examples involve
multipole sources, and their quadrature formulas also contain
derivatives of the solution at the center.
For recent developments of obstacle and partial-balayage
methods in related quadrature-domain settings, see the
Helmholtz quadrature domains of
Kow, Larson, Salo, and Shahgholian \cite{KLSS24},
the Helmholtz partial-balayage results of Gardiner and
Sj\"odin \cite{GardSj}, and the multiphase subharmonic
quadrature domains of Kow, Shahgholian, and Sj\"odin
\cite{KSS26}. The results below develop corresponding
potential characterizations and comparison arguments for
symmetric uniformly elliptic divergence-form operators with
bounded measurable coefficients.

In an example due to T.~Lyons and recorded by Koosis
\cite[pp.~78--79]{Koosis83}, one starts with normalized area
measure on the unit disk and moves the mass of an off-center
subdisk to its center $p.$ By the mean value property on the
subdisk, the resulting measure $\mu$ still satisfies
$\int h \, d\mu = h(0)$ for every integrable harmonic
function $h$ on the unit disk. The subharmonic mean value
inequality fails for $u(z) := \log |z - p|,$ since the atom
at $p$ gives $\int u \, d\mu = -\infty < u(0).$
Brummelhuis and Wiegerinck \cite[Example~2.9]{BW91} give
a related example for the disk algebra.

Hedenmalm, Jakobsson, and Shimorin study weights on the unit
disk for which the weighted average of every bounded harmonic
function equals its value at the origin
\cite[Equation~(1.7)]{HJS02}. For positive real-analytic
weights with subharmonic logarithm, Hedenmalm and Shimorin
show that the harmonic mean value identity determines the
Hele--Shaw domain up to a set of area zero and that this
domain satisfies the corresponding subharmonic inequality
\cite[Theorem~5.3]{HS02}. In our planar example, the change
of variables in Section~\ref{sharpness} gives just such a
one-point harmonic identity, with the Jacobian as the weight,
but the corresponding subharmonic inequality fails.

For positive quadrature data supported at several points,
harmonic quadrature identities need not imply
subharmonic quadrature inequalities even for the
Laplacian. Gustafsson constructed such an example with three
positive point masses \cite[Section~4]{Gus88}. See also
\cite[Example~1.1]{AHT} for a geometric explanation of the
failure of the subharmonic inequality.
Our Theorem~\ref{WDNIS} establishes such a failure
for a single positive point mass by allowing smooth
variable coefficients.
Related uniqueness results for quadrature domains of the
Laplacian were obtained by Shahgholian under conditions on
their modified Schwarz potentials and an internally tangent
ball condition on the domains; see
\cite[Theorem 5.1 and Corollary 5.2]{Shah}.
His setting allows compactly supported distributions as
quadrature data, whereas we use a positive point mass and
deduce nonnegativity of the associated modified Schwarz
$L$-potential from the strong mean value inequality.

Our potential construction and comparison arguments draw
on the classical theory of quadrature domains and Schwarz
potentials; see \cite{S} and \cite[Sections 4.2--4.3]{Sh}.
The uniqueness theorem stated as Theorem 4.7 in Shapiro's
book is due to Epstein and Schiffer \cite{EpSchiff}.
The proof presented by Shapiro is based on an idea of
Sakai and follows a presentation by Gustafsson, as
acknowledged in \cite[p.~30]{Sh}.

Section~\ref{AND} gives the standing assumptions and notation.
Section~\ref{MVTOP} describes the obstacle construction of the
canonical sets and proves the localized mean value inequality.
Section~\ref{ACF} uses the Dirichlet Green function to construct
the weak modified Schwarz $L$-potential of a weak mean value set.
Section~\ref{TNotSP} relates this potential to the signed
mean-value gap and uses an energy comparison to control the
symmetric difference from the canonical set in terms of the
maximal gap. For the variational framework, see \cite{KS}.
Section~\ref{TUT} classifies the open strong mean value sets
and gives conditions under which weak mean value sets are strong.
Section~\ref{sharpness} constructs smooth counterexamples in
every dimension $n \ge 2$ and gives an explicit planar criterion.
Section~\ref{holes} proves the enclosure theorem and uses it to
show that boundary conductivity measurements do not determine
whether a canonical mean value set has a hole when averages are
taken with respect to Lebesgue measure.
The formation, evolution, and extinction of holes
in these canonical families are studied in the
companion paper \cite{ABpartII}.

We benefitted tremendously from the use of AI.
We used AI to assist with developing and checking mathematical
arguments, checking calculations and references, revising
the exposition, creating the figures, and especially 
developing the counterexample in
Section~\refthm{sharpness}\!\!. We take full
responsibility for the content of this paper, including all
mathematical statements, proofs, calculations, and references,
and all AI-generated mathematical suggestions, calculations,
citations, and factual claims used in the manuscript were
independently checked (and frequently simplified) by the authors.

We wish to thank Tim Mesikepp for pointing us in the
direction of quadrature domains and Dave Auckly, Tom Chunders,
and Mike Jensen for helpful discussions.  We thank Jon Tyson
for helping us make better use of AI.

\section{Assumptions, Notation, Definitions, and Background} \label{AND}

In this section we collect our standing coefficient and
auxiliary-domain assumptions,
notation, and references to the definitions introduced
above, together with the Green representation and local
elliptic estimates used later. We begin with the
assumptions on the operator $L.$

We work in $\R^n$ with $n \ge 2.$ Whenever $\Omega$ is used
as an auxiliary domain (also called an ambient domain)
for a Green's function or an obstacle problem, we assume
that it is open, regular for the Dirichlet problem for $L,$
and admits a Dirichlet Green's function. We refer to these
requirements as the \textit{standing auxiliary-domain assumptions}.
For example, bounded $C^{1,1}$ domains satisfy these assumptions;
cf.~\cite[Sections~5--6 and Corollary~9.1]{LSW}.
Boundedness of $\Omega$ is imposed where stated below.

In the definitions of weak and strong mean value sets,
the test functions belong to
$L^1(D)\cap W^{1,2}_{\mathrm{loc}}(D)$ and satisfy the relevant
equation or inequality in $D.$ No extension beyond $D$
and no boundary regularity of $D$ are assumed.
For a fixed $D,$ an auxiliary domain $\Omega$ containing $D$ imposes
no additional conditions on these test functions; they need not extend
to $\Omega.$

We use Einstein summation notation, and we let
\begin{equation}
   Lu := \partial_i(a^{ij}(x)\partial_j u) = \div (A(x) \nabla u) \,.
\label{eq:Ldef}
\end{equation}
We use the (analyst's) convention
$\Delta u = \sum_{i=1}^n \partial_{ii}u.$
For $u\in W^{1,1}_{\mathrm{loc}}(\Omega)$, we interpret $Lu = f$
\emph{distributionally on the open set $\Omega$}:
for every $\varphi\in C_c^\infty(\Omega),$
\begin{equation}
   \int_\Omega (A(x)\nabla u(x))\cdot\nabla\varphi(x)\,dx
 = -\langle f,\varphi\rangle.
\label{eq:Ltestdef}
\end{equation}
When $f$ is a locally integrable function, the pairing on the
right is given by
$$\langle f,\varphi \rangle
 = \int_{\Omega} f(x)\varphi(x) \, dx.$$
If, in addition, $u \in W^{1,2}(\Omega)$ and
$f \in L^2(\Omega),$ then the identity extends by density
to every $\varphi \in W^{1,2}_0(\Omega).$

Equations involving a Dirac mass are understood
distributionally, with smooth compactly supported
test functions. When the Dirac masses cancel upon subtracting two
such equations, as in the proof of Theorem~\refthm{SUR}\!\!, we may use
$W^{1,2}_0(\Omega)$ test functions after verifying the preceding
hypotheses for the difference equation.
Let $D \Subset \Omega$ be open and let $x_0 \in D.$ Deleting
a point $y \in D \setminus \{x_0\}$ leaves Lebesgue integrals
unchanged but makes $G(\cdot,y)$ an admissible $L$-harmonic
test function on $D \setminus \{y\}.$ This can destroy the weak
mean value property at $x_0,$ as Proposition\refthm{pointdefect}and
Remark\refthm{deletepoint}show.

Throughout the paper, unless explicitly qualified by
``almost everywhere,'' a pointwise equality or inequality
stated on a set is understood to hold at every point of
that set, for the representative under consideration.

Finally, we assume that $A$ is bounded and measurable,
symmetric, and uniformly elliptic on the open set under
consideration. More precisely, we assume that there exist
constants $0 < \lambda \le \Lambda < \infty$ such that,
for almost every $x$ in that open set and every
$\eta \in \R^n \setminus \{0\},$
\begin{equation}
   \lambda |\eta|^2 \le a^{ij}(x) \eta_i \eta_j
 \le \Lambda |\eta|^2 \;, \ \text{and} \
                        a^{ij}(x) \equiv a^{ji}(x) \;.
\label{eq:unifellip}
\end{equation}
We refer to these requirements as the
\textit{standing coefficient assumptions}.
We defined weak and strong mean value sets above in
Definitions\refthm{MVsets}and\refthm{SMVsets}\!\!.

For any $f \in L^{\infty}(\Omega),$ when $\Omega$ is bounded, we define
\begin{equation}
(\mathcal{G}f)(x) := \int_{\Omega} G(x,y) f(y) \; dy \;,
\label{eq:Goperator}
\end{equation}
where $G(x,y)$ is the usual Green function for $L$ on $\Omega.$
Fix $p = n+1.$ For $h \in L^p(\Omega),$ let $Th$ be the
unique solution in $W^{1,2}_0(\Omega)$ of $-L(Th) = h.$
The energy estimate and the interior \Holder estimate
of \cite[Chapter~8]{GT} give
$|Th(x)| \le C_x\|h\|_{L^p(\Omega)}$ for the continuous
representative of $Th$ at each $x \in \Omega.$
Fix such an $x.$ By $L^p$ duality there is a unique
$g_x \in L^{p/(p-1)}(\Omega)$ such that
$$Th(x) = \int_\Omega g_x(y)h(y)\,dy
  \qquad \text{for every } h \in L^p(\Omega).$$
For $h \in C_c^\infty(\Omega),$ this is precisely the
duality characterization of the Dirichlet Green function
with pole at $x$ (cf.~\cite[Definition~5.1 and p.~62]{LSW}).
Thus $g_x(y) = G(y,x) = G(x,y)$ for almost every $y,$
where the second equality uses symmetry. Since $\Omega$
is bounded, we may take $h = f$ for every
$f \in L^\infty(\Omega),$ giving
$(Tf)(x) = (\mathcal{G}f)(x)$ at every interior point.
\begin{proposition}[Green representation]\label{greenfacts}
Suppose $\Omega$ is a bounded auxiliary domain satisfying the
standing auxiliary-domain assumptions and $f \in L^\infty(\Omega).$
The function
\begin{equation}\label{eq:psidef}
  \psi_f(x) := -(\mathcal{G}f)(x)
\end{equation}
is the unique weak solution in $W^{1,2}_0(\Omega)$ of
\begin{equation}\label{eq:psieq}
  L\psi_f = f\quad\hbox{in}\quad\Omega.
\end{equation}
At every $x \in \Omega,$ the negative of the integral
in\refeqn{Goperator}gives the value of the locally
\Holder continuous representative of $\psi_f.$
\end{proposition}
We also record the standard local estimates used in the proof of
Lemma~\refthm{localmvi}\!\!\!. We state them only for our
operator $L=\div(A\nabla)$, with bounded measurable,
symmetric coefficients satisfying \eqref{eq:unifellip}.

\begin{lemma}[Local boundedness]\label{localbounds}
Let $B_{2R}(y)$ be a ball on which the standing coefficient
assumptions hold, and write $B_\rho:=B_\rho(y)$.
Suppose that $z\in W^{1,2}(B_{2R})$ is nonnegative and
satisfies $Lz\ge0$ weakly in $B_{2R}$. Then
$$\|z\|_{L^\infty(B_R)}
  \le C R^{-n/2}\|z\|_{L^2(B_{2R})}.$$
Moreover,
$$\|z\|_{L^\infty(B_R)}
  \le C R^{-n}\int_{B_{2R}}z\,dx.$$
Here $C$ depends only on $n$, $\lambda$, and $\Lambda$,
and the $L^\infty$ norms are essential suprema.
\end{lemma}
\noindent
The first estimate is the homogeneous $p=2$ case of
\cite[Theorem 8.17]{GT}. The second is also stated directly
in \cite[Theorem 4.1]{HL}, with $p=1,$ after translation
and rescaling.

\begin{lemma}[Caccioppoli estimates]\label{caccioppoli}
Let $B_{2R}(y)$ be a ball on which the standing coefficient
assumptions hold, and write $B_\rho:=B_\rho(y)$.
If $z\in W^{1,2}(B_{2R})$ is nonnegative and satisfies
$Lz\ge0$ weakly in $B_{2R}$, then
$$\|\nabla z\|_{L^2(B_R)}
  \le C R^{-1}\|z\|_{L^2(B_{2R})}.$$
If $u\in W^{1,2}(B_{2R})$ satisfies $Lu=f$ weakly
in $B_{2R}$, with $f\in L^\infty(B_{2R})$, then
$$\|\nabla u\|_{L^2(B_R)}
  \le C\left(
       R^{-1}\|u\|_{L^2(B_{2R})}
       +R\|f\|_{L^2(B_{2R})}
       \right).$$
In both estimates, $C$ depends only on
$n$, $\lambda$, and $\Lambda$.
\end{lemma}
For these standard cutoff energy estimates, see
\cite[Theorem 2.1, Theorem 3.36, and Remark 3.40]{AMC}.

We use standard notation for balls and Sobolev spaces;
see Gilbarg and Trudinger \cite{GT}.
The symbols $\operatorname{int} E,$ $\overline{E},$
$D\triangle E,$ $\bigrchi_E,$ and $\delta_y$
denote interior, closure, symmetric difference,
characteristic function, and unit Dirac mass, respectively.
For an open set $U,$ we write $E \Subset U$ when
$\overline{E}$ is a compact subset of $U,$
and use $\subset$ for inclusion allowing equality.
We set $\omega_n := |B_1(0)|.$
Table~\ref{tab:notation} collects the notation for the
operator, Green functions, and mean value sets,
including objects introduced in the following sections.

\begin{table}[H]
\centering
\caption{Standing assumptions and frequently used notation}
\label{tab:notation}
\small
\renewcommand{\arraystretch}{1.25}
\begin{tabular}{@{}p{0.15\textwidth}p{0.62\textwidth}p{0.15\textwidth}@{}}
\hline
\multicolumn{3}{@{}l}{\textbf{Standing assumptions}} \\[3pt]
\textbf{Applies to} & \textbf{Requirements} & \textbf{Reference} \\
\hline
Coefficients
  & $A$ is bounded and measurable, symmetric, and uniformly
    elliptic with constants $0<\lambda\le\Lambda<\infty.$
  & \eqref{eq:unifellip} \\[4pt]
Auxiliary domains
  & $\Omega$ is open, regular for the Dirichlet problem for $L,$
    and admits a Dirichlet Green function. Boundedness is imposed
    separately when needed.
  & Section~\ref{AND} \\[4pt]
\hline
\multicolumn{3}{@{}l}{\textbf{Frequently used notation}} \\[3pt]
\textbf{Notation} & \textbf{Meaning} & \textbf{Reference} \\
\hline
$L, \, A$
  & $Lu=\operatorname{div}(A\nabla u),$ where
    $A(x) = (a^{ij}(x))$ is the coefficient matrix.
  & \eqref{eq:Ldef} \\[4pt]
$\lambda, \, \Lambda$
  & Lower and upper uniform ellipticity constants.
  & \eqref{eq:unifellip} \\[4pt]
$\Omega$
  & Auxiliary domain for Green functions and obstacle problems.
  & Section~\ref{AND} \\[4pt]
$G(x,y)$
  & Positive Dirichlet Green function, with
    $L_x G(x,y) = -\delta_y.$
  & Section~\ref{AND} \\[4pt]
$\mathcal{G}$
  & Green operator:
    $\mathcal{G}f(x) := \int_\Omega G(x,y)f(y) \, dy.$
  & \eqref{eq:Goperator} \\[4pt]
$\psi_f$
  & Zero-boundary solution of $L\psi_f=f$ in $\Omega$:
    $\psi_f := -\mathcal{G}f.$
  & \eqref{eq:psidef} \\[4pt]
$D_r(x_0)$
  & Canonical mean value set centered at $x_0,$ with
    $|D_r(x_0)| = r^n$; $r$ is a volume parameter.
  & Theorem~\ref{KeyMVT} \\[4pt]
$u_r$
  & Solution of the canonical obstacle problem.
  & \eqref{eq:wrobprob} \\[4pt]
$W_r$
  & Canonical height: $W_r:=G(\cdot,x_0)-u_r.$
  & \eqref{eq:hfdef} \\[4pt]
$W_D$
  & Weak modified Schwarz $L$-potential associated
    with a weak mean value set $D$ at $x_0.$
  & \eqref{eq:potdef} \\
$\Theta_D(u)$
  & Signed mean-value gap:
    $\Theta_D(u) := u(x_0) - \dfrac{1}{|D|}\int_D u(x) \, dx.$
  & \eqref{eq:mvgap} \\
$\Upsilon_{D,E}$
  & Canonical comparison:
    $\Upsilon_{D,E}:=W_D-W_E,$ where $E$ is the canonical
    set with the same center and volume as $D.$
  & \eqref{eq:can} \\[4pt]
\hline
\end{tabular}
\end{table}

\section{The Mean Value Theorem via the Obstacle Problem}
\label{MVTOP}

In the Fermi Lectures \cite{C2}, Caffarelli presented
an obstacle problem approach to the mean value theorem
and proposed its extension to uniformly elliptic
divergence form operators.
The key observation is that differences of suitably
ordered obstacle solutions give nonnegative, compactly
supported test functions whose equations compare
averages over the corresponding noncontact sets.
Blank and Hao \cite{BH1} established this construction
for operators with bounded measurable coefficients.
Although the presentation in \cite{BH1} is geared toward
$n > 2,$ the results used below hold for all $n \ge 2.$
We use the bounded-domain Dirichlet Green function construction
of Benson, Blank, and LeCrone \cite[Section 4]{BBL}.
Subtracting the $L$-harmonic extension of the common boundary
data from both the fundamental solution and the obstacle minimizer
preserves their difference, and hence the noncontact set.  By
doing it this way, only the proof of Theorem\refthm{KeyMVT}\!\!(2)
is different for $n = 2.$ The estimates for (2) in this case are
given in Remark~\refthm{GFbeatsFS}\!\!.

The formulation below includes a modest but useful
strengthening of the mean value inequalities stated in
\cite[Theorem 6.3]{BH1} and \cite[Theorem 1.1]{BBL}:
the test function need only be defined on the averaging
set itself.  Specifically, in part~(3) we require only
$$v \in L^1(D_s(x_0))\cap W^{1,2}_{\mathrm{loc}}(D_s(x_0)),
  \ \ \ \ Lv \ge 0\ \hbox{weakly in }D_s(x_0),$$
with no extension of $v$ as a subsolution to the
auxiliary domain $\Omega.$  This localized assertion,
including the endpoint $s = r_0,$ is proved in
Lemma~\refthm{localmvi}below.  This refinement is needed
to conclude that the canonical sets are strong mean value
sets under our definition, which tests subsolutions defined
only on the set itself.

\begin{theorem} [Mean value theorem for divergence form
                     elliptic PDE]  \label{KeyMVT}
Let $L$ be defined as in Section\refthm{AND}in the open connected
bounded domain $\Omega.$ For any $x_0 \in \Omega,$ define $D_r(x_0)$ to
be the noncontact set for the solution in $W^{1,2}_0(\Omega)$ to the
obstacle problem:
\begin{equation}
      \begin{array}{rll}
             L u_r &\!\!\! = -r^{-n}
               \bigrchi_{\left\lbrace u_r < G(\cdot, x_0) \right\rbrace} \ \ \
                 & \text{in} \ \Omega \\
              u_r &\!\!\! \le G(\cdot, x_0) \ \ \ & \text{in} \ \Omega
      \end{array}
\label{eq:wrobprob}
\end{equation}
and let $r_0 = r_0(\Omega,x_0)$ denote the supremum of those
$r > 0$ for which $D_r(x_0) \Subset \Omega.$ Then
$r_0 > 0,$ and for $r \in (0, r_0]$ the sets $\{ D_r(x_0) \}$ have
the following properties:
\begin{enumerate}
  \item If $0 < r < s \le r_0,$ then $D_r(x_0) \subset D_s(x_0).$
  \item We have the inclusions: $B_{cr}(x_0) \subset D_r(x_0)
           \subset B_{Cr}(x_0),$
           with $c,$ $C$ depending only on $n,$ $\lambda,$
           and $\Lambda.$
\item For any $r$ and $s$ satisfying $0 < r < s \le r_0,$
      and for any
      $$v \in L^1(D_s(x_0))
          \cap W^{1,2}_{\mathrm{loc}}(D_s(x_0))$$
      satisfying $Lv \ge 0$ weakly in $D_s(x_0),$ we have
      \begin{equation}
             v(x_0)
             \le \frac{1}{|D_r(x_0)|}\int_{D_r(x_0)} v
             \le \frac{1}{|D_s(x_0)|}\int_{D_s(x_0)} v \;.
      \label{eq:MVTres}
      \end{equation}
      If $Lv \le 0$ weakly in $D_s(x_0),$ then the
      inequalities in Equation~\refeqn{MVTres}are reversed.
      Consequently, equality holds throughout when
      $Lv = 0$ weakly in $D_s(x_0).$
      In particular, the sets $D_r(x_0)$ are strong
      mean value sets.
  \item If $\Omega_1$ and $\Omega_2$ are bounded auxiliary
        domains containing $x_0$ for the same operator $L,$
        then they give the same set $D_s(x_0)$ whenever
        $0 < s \le \min\{r_0(\Omega_1,x_0),r_0(\Omega_2,x_0)\}.$
  \item Finally, we have $|D_r(x_0)| = r^n$ and $| \partial D_r(x_0) | = 0.$
\end{enumerate}
The set $D_r(x_0)$ is called the
\textit{canonical mean value set} for $L$ at $x_0$
with \textit{canonical parameter} $r$ and
volume $r^n.$
\end{theorem}
\begin{figure}[H]
\centering
\begin{minipage}[c]{0.50\textwidth}
  \centering
  \includegraphics[width=\linewidth]{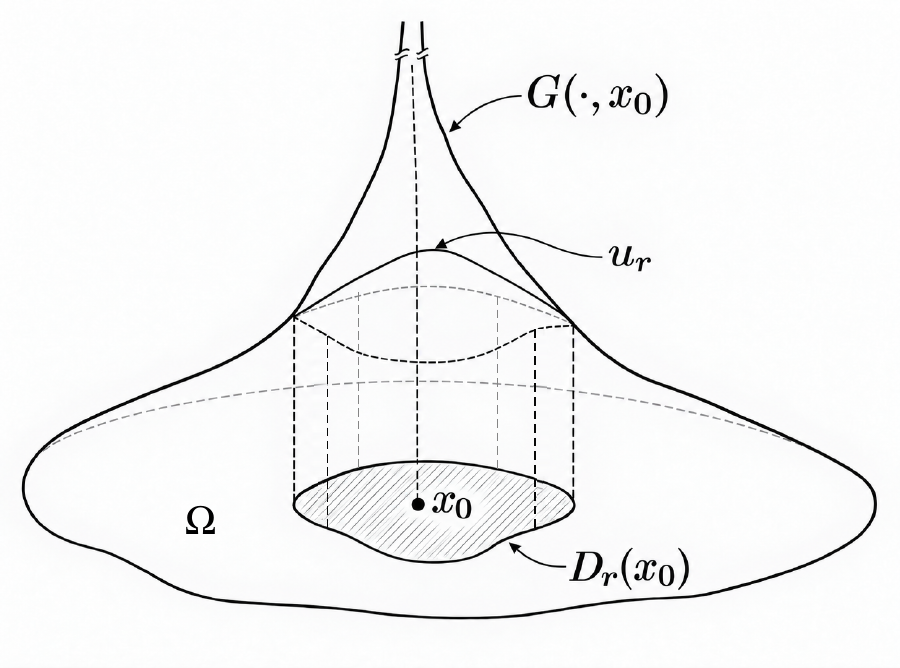}
\end{minipage}\hspace{0.04\textwidth}%
\begin{minipage}[c]{0.36\textwidth}
  \setlength{\captionindent}{0pt}
  \setlength{\abovecaptionskip}{0pt}
  \caption{Schematic illustration of the canonical obstacle
  construction in Theorem~\ref{KeyMVT}.
  The height function $W_r=G(\cdot,x_0)-u_r$ measures
  the vertical gap between the graphs and vanishes
  on $\Omega\setminus D_r(x_0)$.}
  \label{fig:green-obstacle}
\end{minipage}
\end{figure}

\begin{remark}[Existence and discussion of solutions to
        Equation\refeqn{wrobprob}\!\!]    \label{existence}
The existence of solutions to Equation\refeqn{wrobprob}follows from
the theorems within \cite{BH1}. In terms of
understanding this equation, one can describe the solution
(heuristically) as being a membrane that was inflated from below
the Green's function; see Figure~\ref{fig:green-obstacle}.
The parameter $r$ is related to the pressure of the inflation.
This description is only heuristic, but suggests correctly that
the problem can be reformulated as a convex minimization problem
within the calculus of variations.
\end{remark}

\begin{remark}[Green's functions and the ambient domain]
\label{GFbeatsFS}
Caffarelli and Blank--Hao use a fundamental solution on
$\R^n$ in \cite{C2,BH1}. Following \cite{BBL}, we instead
use the Dirichlet Green's function on an auxiliary domain
$\Omega$, as in \eqref{eq:wrobprob}.
The independence of the canonical sets from the choice
of sufficiently large ambient domain, asserted in
Theorem~\ref{KeyMVT}(4), is justified by the
harmonic-correction argument and its extension to
$r = r_0$ in Lemma~\ref{mvprops} below.
When considering the family $D_r(x_0)$ for all $r > 0,$
we fix the coefficient field $A$ on all of $\R^n,$
as in Theorem~\refthm{openreps}\!\!.
Indeed, for any prescribed $r > 0,$ the canonical set can be
constructed in a sufficiently large ball centered at $x_0$
and is compactly contained in that ball. For $n > 2,$ this
follows from \cite[Theorem~6.3]{BH1}.

In dimension two, once the ball inclusions in part~(2) of
Theorem~\refthm{KeyMVT}are established, the arguments for the
remaining parts apply without further modification. To prove
the inclusions, scale to $r = 1$ and take a disk $B_M$ centered
at the pole $x_0.$ On $B_{3M/4} \setminus B_{M/4},$ the height
is bounded above by the Dirichlet Green function, which by
\cite[Theorem~7.1]{LSW} satisfies
$G_{B_M}(x,x_0) \le C \log(M/|x - x_0|) \le C \log 4.$
Positivity at any point of $\partial B_{M/2}$ would force the
height to reach $cM^2$ on the annulus by the nondegeneracy
estimate of Blank and Hao \cite[Theorem~3.9]{BH1}, so for large
$M$ the maximum principle confines the noncontact set to $B_{M/2}$.
With $M$ now fixed, the logarithmic lower bound near the pole
dominates the local upper bound for the minimizer, giving the
inner ball as in \cite[Lemmas~4.4 and~6.2]{BH1}. Rescaling gives
$B_{cr}(x_0) \subset D_r(x_0) \subset B_{Mr/2}(x_0) \Subset B_{Mr}(x_0).$
Lemma~\refthm{mvprops}then shows that the resulting
canonical set and height are independent of the
particular sufficiently large auxiliary domain.
\end{remark}
\begin{remark}[The height function]   \label{height}
Because it will come up later, we note that the ``height function''
\begin{equation}
\label{eq:hfdef}
W_r(x) := G(x,x_0) - u_r(x)
\end{equation}
satisfies the following:
\begin{equation}
      \begin{array}{rll}
          LW_r &\!\!\! = r^{-n} \bigrchi_{\left\lbrace W_r > 0 \right\rbrace}
                        - \delta_{x_0} \ \ \ & \text{in} \ \Omega \\
           W_r &\!\!\! \ge 0 \ \ \ & \text{in} \ \Omega \\
           W_r &\!\!\! = 0 \ \ \ \ &\text{on} \ \partial \Omega \;.
      \end{array}
\label{eq:height}
\end{equation}
\end{remark}

\begin{lemma}[Volume, continuity, and ambient domains]
\label{mvprops}
Use the notation and hypotheses of
Theorem~\ref{KeyMVT}.
For every $0 < r \le r_0,$
$$|D_r(x_0)| = r^n \ \ \ \ \text{(``Mass Balance'').}$$
As $r \uparrow r_0,$ the obstacle solutions $u_r$
converge to $u_{r_0}$ in $W^{1,2}_0(\Omega)$ and
locally uniformly in $\Omega,$ and
$$D_{r_0}(x_0)
  = \bigcup_{0 < r < r_0} D_r(x_0).$$
The height functions and canonical sets are independent
of the choice of sufficiently large auxiliary domain,
including at $r = r_0.$ Moreover,
$$|\partial D_{r_0}(x_0)| = 0.$$
\end{lemma}

\begin{pf}
\textit{Mass balance.}
\nopagebreak[4]

Let $D_r(x_0) \Subset \Omega$ and let $\eta \in C^{\infty}_c(\Omega)$
satisfy $\eta \equiv 1$ on a neighborhood of $\pclosure{D_r(x_0)}.$
Since $\eta$ is constant on a neighborhood of
$\overline{D_r(x_0)}$, its gradient is supported away from
$\overline{D_r(x_0)}$ and therefore also away from $x_0$.  On that
support $W_r$ vanishes,
and $W_r$ is locally in $W^{1,2}$ away from $x_0$.  Consequently,
$\nabla W_r = 0$ almost everywhere on $\operatorname{supp}\nabla\eta$,
and the weak formulation gives
$$\langle LW_r,\eta \rangle
 = -\int_{\Omega} A\nabla W_r\cdot\nabla\eta \, dx = 0.$$
On the other hand, using Equation\refeqn{height}we have
$$\langle L W_r, \eta \rangle =
       \langle r^{-n}\bigrchi_{\left\lbrace W_r > 0 \right\rbrace}
                        - \delta_{x_0}, \eta \rangle
 = r^{-n}|\{ W_r > 0 \}| - 1 \;,$$
so $r^n = |D_r(x_0)| \le |\Omega|$ whenever
$D_r(x_0) \Subset \Omega.$ Taking the supremum over these
$r$ gives $r_0 \le |\Omega|^{1/n} < \infty.$

\medskip
\noindent
\textit{Continuity at the endpoint.}
\nopagebreak[4]

We sketch the parameter-continuity argument, adapting
\cite[Lemma 2.6 and Remark 2.7]{AryB} to the present
Green-function formulation; see also \cite[Lemma 4.5]{BBL}.
The argument uses the variational formulation and interior
estimates developed in \cite{BH1}.  Assume $\Omega$ is a fixed bounded
ambient domain.  If $r,s$ stay in a compact subinterval of
$(0,\infty)$, we may use each solution as a competitor
for the other, since both belong to $W^{1,2}_0(\Omega)$
and satisfy the same obstacle constraint
$u \le G(\cdot,x_0).$ The two variational inequalities give
$$\begin{aligned}
  \int_\Omega A\nabla u_r \cdot \nabla(u_s-u_r)\,dx
    &\ge r^{-n}\int_\Omega(u_s-u_r)\,dx,\\
  \int_\Omega A\nabla u_s \cdot \nabla(u_r-u_s)\,dx
    &\ge s^{-n}\int_\Omega(u_r-u_s)\,dx.
\end{aligned}$$
Combining these inequalities, writing $w := u_r-u_s,$ and
using the Cauchy--Schwarz inequality, we obtain
$$\int_\Omega A\nabla w\cdot\nabla w\,dx
  \le (r^{-n}-s^{-n})\int_\Omega w\,dx
  \le |r^{-n}-s^{-n}|\,
          |\Omega|^{1/2}\|w\|_{L^2(\Omega)}.$$
Combining that estimate with ellipticity and the
Poincar\'e inequality gives us
$$\begin{aligned}
  \lambda\|\nabla w\|_{L^2(\Omega)}^2
    &\le \int_\Omega A\nabla w\cdot\nabla w\,dx\\
    &\le |r^{-n}-s^{-n}|\,
          |\Omega|^{1/2}\|w\|_{L^2(\Omega)}\\
    &\le C|\Omega|^{1/2}|r^{-n}-s^{-n}|\,
          \|\nabla w\|_{L^2(\Omega)}.
\end{aligned}$$
If $\|\nabla w\|_{L^2(\Omega)} = 0,$ then $w = 0.$
Otherwise, dividing by this norm and applying the
Poincar\'e inequality once more gives
$$\|u_r-u_s\|_{W^{1,2}_0(\Omega)}
  \le C\,|r^{-n}-s^{-n}|,$$
where $C$ depends only on $\lambda$ and $\Omega,$
and is independent of $r$ and $s.$
For $r$ in a compact interval $I \subset (0,\infty)$ containing
$r_0$, the right-hand sides
$-r^{-n}\bigrchi_{\left\lbrace D_r(x_0) \right\rbrace}$ are uniformly
bounded in $L^\infty(\Omega)$.
The preceding $W^{1,2}_0(\Omega)$ estimate gives a uniform
$W^{1,2}$ bound for $u_r$ on $I$.  Consequently, the interior
De Giorgi--Nash--Moser estimates give, for every $K\Subset\Omega$,
constants $\alpha \in (0,1)$ and $C_K$ independent of $r \in I$
such that
$$ \|u_r\|_{C^\alpha(K)} \le C_K .$$
The $W^{1,2}$ convergence established above, together with this
uniform H\"older bound, implies that $u_r\to u_{r_0}$ uniformly
on every such $K$ as $r\to r_0$.  Thus $W_r = G(\cdot,x_0)-u_r$
converges locally uniformly away from $x_0$;
the Green-function singularity ensures that $x_0$ belongs to every
noncontact set under consideration.

For $0 < r < s,$ the variational comparison in the proof of
\cite[Lemma~6.1]{BH1} gives $u_r \ge u_s$ and hence $W_r \le W_s.$
Combining this ordering with the convergence as $r \uparrow r_0,$
we obtain
\begin{equation}
 D_{r_0}(x_0) = \{W_{r_0} > 0\}
 = \bigcup_{0 < r < r_0}\{W_r > 0\}
 = \bigcup_{0 < r < r_0}D_r(x_0).
\label{eq:union}
\end{equation}
In particular, the volume identity proved above for
$r < r_0$ extends to $r = r_0$ by continuity of measure
from below.

\medskip
\noindent
\textit{Independence of the ambient domain.}
\nopagebreak[4]

Let $D_t^{(i)}$ denote the set constructed in $\Omega_i,$
with the same operator $L$ and center $x_0$ in both domains.
Fix $0 < s \le \min\{r_0(\Omega_1,x_0),r_0(\Omega_2,x_0)\}.$
For every $0 < t < s,$ we have $D_t^{(i)} \Subset \Omega_i$
for $i = 1,2.$ Choose a ball $B$ with
$\Omega_1 \cup \Omega_2 \Subset B.$ If $A$ is given only
on this union, extend it by $\lambda I$ elsewhere in $B.$

For $U = B$ or $\Omega_i,$ write $G_U$ for the Dirichlet Green
function with pole $x_0.$ Cutting off $G_B$ near $x_0$ gives a
$W^{1,2}(\Omega_i)$ function agreeing with $G_B$ near
$\partial\Omega_i.$ Let $h_i \in W^{1,2}(\Omega_i)$ be its
$L$-harmonic replacement. To identify $h_i,$ replace the pole
in both Green problems by the same smooth approximate identity
and denote the resulting solutions by $g_{U,\epsilon}.$
Green representation and continuity off the diagonal give
$g_{U,\epsilon} \to G_U$ locally uniformly away from $x_0$
\cite[Sections~6--7]{LSW}. In particular, $g_{B,\epsilon} \to G_B$
uniformly on $\partial\Omega_i.$ Since
$g_{B,\epsilon} - g_{\Omega_i,\epsilon}$ is $L$-harmonic in
$\Omega_i$ with boundary values $g_{B,\epsilon},$ the maximum
principle gives uniform convergence of this difference to $h_i.$
Thus $G_B - G_{\Omega_i} = h_i$ away from $x_0.$

As in \cite[Lemma~4.2]{BBL}, add $h_i$ to the obstacle solution
in $\Omega_i$ and extend the sum by $G_B$ to the rest of $B.$
Since $D_t^{(i)} \Subset \Omega_i,$ the sum already equals $G_B$
near $\partial\Omega_i.$ The resulting function solves the
obstacle problem in $B,$ and its height is the zero extension
of the height in $\Omega_i.$ Uniqueness in $B$ therefore gives
the same height for $i = 1,2$ and hence $D_t^{(1)} = D_t^{(2)}.$
The argument giving Equation\refeqn{union}applies at $s$
in each domain, so
$$D_s^{(1)} = \bigcup_{0 < t < s} D_t^{(1)}
            = \bigcup_{0 < t < s} D_t^{(2)} = D_s^{(2)}.$$
Local convergence away from $x_0$ also gives equality of
the two heights on $\Omega_1 \cap \Omega_2$ away from the
pole at parameter $s.$
For $r < r_0,$ the measure-zero free-boundary conclusion
$$|\partial D_r(x_0)| = 0$$
follows from \cite[Corollary 3.10]{BH1}.
For $r = r_0,$ compare with a sufficiently large auxiliary
ball in which the same canonical set is compactly
contained. Applying \cite[Corollary 3.10]{BH1} in that
larger domain gives
$$|\partial D_{r_0}(x_0)| = 0.$$
\end{pf}

The next lemma proves the localized mean value assertion
included in Theorem~\refthm{KeyMVT}\!\!(3).
This localization is needed to conclude that the canonical
sets are strong mean value sets in the sense of
Definition~\refthm{SMVsets}\!\!: the test functions in that
definition need only be integrable subsolutions on the
set itself, with no extension beyond its boundary.
It is also used in the proof of
Theorem~\refthm{openreps}\!\!, where a subsolution on a
given set is restricted to a canonical set contained in it.

If $D_t(x_0) \Subset D_s(x_0)$ for every $0 < t < s \le r_0,$
we could apply the ambient-domain mean value theorem on a smooth
domain $\Omega_t$ with $D_t(x_0) \Subset \Omega_t \Subset D_s(x_0)$
and then let $t \uparrow s.$
Such compact containment holds under stronger regularity assumptions
(see \cite{ArmB}), while the cutoff argument below requires only
our standing assumptions on the coefficients.

\begin{lemma}[Localization of the mean value inequality]
\label{localmvi}
Using the notation and hypotheses of
Theorem~\refthm{KeyMVT}\!\!, we
let $0 < r < s \le r_0,$ and suppose that
$$v \in L^1(D_s(x_0))
       \cap W^{1,2}_{\mathrm{loc}}(D_s(x_0))$$
satisfies $Lv \ge 0$ weakly in $D_s(x_0).$ Then
$$v(x_0)
  \le \frac{1}{|D_r(x_0)|}
      \int_{D_r(x_0)} v(x) \, dx
  \le \frac{1}{|D_s(x_0)|}
      \int_{D_s(x_0)} v(x) \, dx.$$
For supersolutions the inequalities are reversed,
and for solutions both are equalities.
\end{lemma}

\begin{pf}
Fix $0 < r < s < r_0$, write $D := D_s(x_0)$, and set
$$F := W_s - W_r = u_r - u_s, \ \ \ f :=
    s^{-n} \bigrchi_{\left\lbrace D \right\rbrace} -
    r^{-n} \bigrchi_{\left\lbrace D_r(x_0) \right\rbrace}.$$
Then $F \ge 0$, $F \in W^{1,2}_0(\Omega)$, and $LF = f$.
The volume identity in Lemma~\refthm{mvprops}gives
$$\int_D f \, dx = 0.$$
Next, let $v \in L^1(D)\cap W^{1,2}_{\mathrm{loc}}(D)$ satisfy
$Lv \ge 0$ weakly in $D.$ For $M > 0$, put
$$v_M := \max\{v,-M\},\ \ \ z := v_M+M.$$
By the standard truncation property for weak subsolutions,
$v_M$ is again a weak subsolution; equivalently, this
follows by applying the weak inequality to smooth convex
nondecreasing approximations of
$t \mapsto \max\{t,-M\}.$
Hence $z$ is a nonnegative integrable weak subsolution
on $D.$

Write $d(x) := \operatorname{dist}(x,D^c)$ and
$E_\delta := \{x \in D:\delta < d(x) < 2\delta\}$.
Choose compactly supported Lipschitz cutoffs
$\eta_\delta$ on $D$ satisfying
$$0 \le \eta_\delta \le 1,\ \ \
  \eta_\delta = 0\ \hbox{on }\{d \le \delta\},\ \ \
  \eta_\delta = 1\ \hbox{on }\{d \ge 2\delta\},$$
and
$$|\nabla\eta_\delta| \le C\delta^{-1}
  \bigrchi_{\left\lbrace E_\delta \right\rbrace}
  \ \hbox{almost everywhere}.$$

The standard quadratic growth estimate for solutions of the
obstacle problem (see \cite[Theorem 3.2]{BH1} for our
situation) gives
$0 \le F \le W_s \le Cd^2$ near $\partial D.$
Cover $E_\delta$ by balls $B_j$ of radius $\delta/8,$
centered in $E_\delta$, whose concentric enlargements
$4B_j$ have bounded overlap. These enlarged balls lie
in $\{0 < d < 3\delta\}\subset D$.
On $4B_j$, the quadratic bound gives
\begin{equation}\label{eq:local-F-sup}
  \|F\|_{L^\infty(4B_j)}\le C\delta^2.
\end{equation}
Since $LF=f$ and $f$ is bounded, this bound and the
inhomogeneous estimate in Lemma~\ref{caccioppoli} give
\begin{equation}\label{eq:local-F-gradient}
  \|\nabla F\|_{L^2(B_j)}\le C\delta^{n/2+1}.
\end{equation}

For $z$, Lemma~\ref{localbounds} applied to
$2B_j\subset 4B_j$ gives
\begin{equation}\label{eq:local-z-sup}
  \|z\|_{L^\infty(2B_j)}
  \le C\delta^{-n}\int_{4B_j}z\,dx.
\end{equation}
The subsolution estimate in Lemma~\ref{caccioppoli},
together with $|2B_j|^{1/2}\le C\delta^{n/2}$
and \eqref{eq:local-z-sup}, then gives
\begin{equation}\label{eq:local-z-gradient}
  \|\nabla z\|_{L^2(B_j)}
  \le C\delta^{-n/2-1}\int_{4B_j}z\,dx.
\end{equation}

Consequently, the cutoff bound, Cauchy--Schwarz, and
bounded overlap, together with $|B_j|^{1/2}\le C\delta^{n/2}$,
yield
\begin{equation}\label{eq:local-cutoff-error}
\begin{aligned}
&\int_D \bigl(F|\nabla z|+z|\nabla F|\bigr)
              |\nabla\eta_\delta|\,dx\\
&\quad\le C\delta^{-1}\sum_j |B_j|^{1/2}
  \Bigl(
    \|F\|_{L^\infty(4B_j)}\|\nabla z\|_{L^2(B_j)}
    +\|z\|_{L^\infty(2B_j)}\|\nabla F\|_{L^2(B_j)}
  \Bigr)\\
&\quad\le C\sum_j\int_{4B_j}z\,dx\\
&\quad\le C\int_{\{0<d<3\delta\}}z\,dx
  \longrightarrow 0
\end{aligned}
\end{equation}
as $\delta\downarrow0$, by dominated convergence for
fixed $M$. The constants are independent of $\delta$ and $M.$

Since both $F$ and $\eta_{\delta}$ are nonnegative, we can use
$F\eta_{\delta}$ as a test function in $Lz \ge 0$ to get
$$\int_D A \nabla z \cdot \nabla (F \eta_{\delta}) \, dx \le 0,$$
and therefore by the symmetry of $A$ we get
\begin{equation}
\int_D F A \nabla z \cdot \nabla \eta_{\delta} \, dx \le 
  - \int_D \eta_{\delta} A \nabla F \cdot \nabla z \, dx \;.
\label{eq:Lzge0}
\end{equation}
Using $z \eta_{\delta}$ as test function in the equality $LF = f,$
expanding the gradient of the product, and applying
Equation\refeqn{Lzge0}gives 
\begin{alignat*}{2}
\int_D f z \eta_{\delta} \, dx
  &= - \int_D \eta_{\delta} A \nabla F \cdot \nabla z \, dx
     - \int_D z A \nabla F \cdot \nabla \eta_{\delta} \, dx \\
  &\ge \int_D F A \nabla z \cdot \nabla \eta_{\delta} \, dx
     - \int_D z A \nabla F \cdot \nabla \eta_{\delta} \, dx \\
  &= \int_D (FA \nabla z - zA \nabla F) \cdot \nabla \eta_{\delta} \, dx \,.
\end{alignat*}
Since $A$ is bounded, Equation\refeqn{local-cutoff-error}implies
that the right-hand side converges to zero as
$\delta \downarrow 0.$ On the left-hand side,
$\eta_\delta \to 1$ pointwise in $D,$ and
$|z\eta_\delta f|\le\|f\|_{L^\infty(D)}z\in L^1(D).$
Dominated convergence therefore gives
$\int_D zf \, dx \ge 0.$

Next, $\int_D f \, dx = 0,$ so we have
$\int_D v_Mf\,dx \ge 0.$ Letting $M \to \infty$ and
using $|v_M| \le |v|$, we obtain
$$0 \le \int_D vf \, dx
 = \frac{1}{|D_s(x_0)|}\int_{D_s(x_0)}v \, dx
    -\frac{1}{|D_r(x_0)|}\int_{D_r(x_0)}v \, dx.$$
Choose a ball $B$ centered at $x_0$ with $B \Subset D.$
For all sufficiently small $t > 0,$ Lemma~\refthm{mvprops}identifies
$D_t(x_0)$ with the canonical set constructed in $B.$
Applying \cite[Corollary~6.4]{BH1} in $B$ to the supersolution
$-v$ and changing signs gives
$$v(x_0) = \lim_{t \downarrow 0}
  \frac{1}{|D_t(x_0)|}\int_{D_t(x_0)} v \, dx.$$
The monotonicity just proved now gives both inequalities for $s < r_0.$

For $s = r_0$, apply the result to the restriction of
$v$ to $D_t(x_0),$ where $r < t < r_0,$ and let
$t \uparrow r_0.$ Equation\refeqn{union}\!\!,
the volume identity, and dominated convergence give
the endpoint inequalities. Applying the argument to
$-v$ gives the reversed inequalities for supersolutions
and hence equalities for solutions.
\end{pf}

\section{Quadrature Domains and Weak Modified Schwarz $L$-Potentials}
\label{ACF}

Throughout this section, let $\Omega$ be a bounded auxiliary
domain satisfying the standing auxiliary-domain assumptions.
A weak mean value set $D$ for $L$ at $x_0$ is a
quadrature domain with the single positive point-mass
datum $\mu := |D|\delta_{x_0}:$
$$\int_D u(x) \, dx
  = \int_D u(x) \, d\mu(x)
  = |D|u(x_0)$$
for every $u \in L^1(D)\cap W^{1,2}_{\mathrm{loc}}(D)$
satisfying $Lu = 0$ weakly in $D.$

For the Laplacian, Shapiro proves the connection between
quadrature domains and modified Schwarz potentials in
\cite[Theorem 4.1]{Sh}.
Our construction adapts the preliminary potential argument
in \cite[Section 4.2, p.~26]{Sh}, using the Dirichlet Green
function of $L$ in place of the fundamental solution of
the Laplacian.
This gives the following analogue for weak mean value sets
and the operator $L$:

\begin{theorem}[Weak mean value sets yield weak modified
Schwarz $L$-potentials]
\label{MVSYWMSP}
Let $D \Subset \Omega$ be a weak mean value
set for $L$ at $x_0 \in D.$ Then there exists a unique
function $W_D$ in the class
$$W_D - G(\cdot,x_0) \in W^{1,2}_0(\Omega)$$
such that
\begin{equation}
\begin{aligned}
LW_D &= \frac{1}{|D|}\bigrchi_{D} - \delta_{x_0}
  &&\text{in } \Omega,\\
W_D &= 0
  &&\text{on } \Omega \setminus D.
\end{aligned}
\label{eq:OPHF}
\end{equation}
The differential equation is understood
distributionally.
We use the continuous representative of $W_D - G(\cdot,x_0)$
to specify $W_D$ pointwise on $\Omega \setminus \{x_0\}.$
\end{theorem}

\begin{pf}
Let $G$ be the positive Dirichlet Green function
for $L$ on $\Omega,$ normalized by
$$LG(\cdot,y) = -\delta_y.$$
Define
\begin{equation}
W_D(z) := G(z,x_0)
  - \frac{1}{|D|}\int_D G(z,x) \, dx.
\label{eq:potdef}
\end{equation}
Fix $z \in \Omega \setminus D.$ By symmetry of the
Green function, $x \mapsto G(z,x)$ is $L$-harmonic
in $D.$ Moreover, $G(z,\cdot) \in L^1(D),$ even
when $z \in \partial D,$ because the Green-function
singularity is locally integrable and
$D \Subset \Omega.$ Thus the weak mean
value property gives
$$\frac{1}{|D|}\int_D G(z,x) \, dx = G(z,x_0).$$
Substituting this identity into
Equation~\refeqn{potdef}shows that
$$W_D(z) = 0
  \qquad \text{for every } z \in \Omega \setminus D.$$

Set $f := \bigrchi_D/|D|.$ By Proposition~\refthm{greenfacts}\!\!,
$$W_D - G(\cdot,x_0) = \psi_f \in W^{1,2}_0(\Omega).$$
Since $L\psi_f = f$ and $LG(\cdot,x_0) = -\delta_{x_0},$
the distributional equation in\refeqn{OPHF}follows.

To prove uniqueness, suppose that
$\widetilde{W}_D$ is another function in the stated
class satisfying Equation~\eqref{eq:OPHF}.
Then
$$q := W_D - \widetilde{W}_D
  \in W^{1,2}_0(\Omega), \qquad Lq = 0.$$
Testing with $q$ and using uniform ellipticity gives
$$0 = \int_{\Omega} A\nabla q \cdot \nabla q \, dx
  \ge \lambda\int_{\Omega} |\nabla q|^2 \, dx.$$
Thus $q = 0$ by the Poincar\'e inequality.
\end{pf}

We call $W_D$ the \textit{weak modified Schwarz
$L$-potential} of $D$ at $x_0.$
Nonnegativity is not part of this definition,
and no additional condition requiring
$\nabla W_D = 0$ on $\partial D$ is imposed.

\section{Comparison of Schwarz Potentials}
\label{TNotSP}

For a bounded open set $D$ containing $x_0,$ and for
$u\in L^1(D)$ with a specified pointwise value at $x_0,$
we define the \textit{signed mean-value gap} of $u$ on $D$ by
\begin{equation}
\Theta_D(u)
  := u(x_0)-\frac{1}{|D|}\int_D u(x)\,dx.
\label{eq:mvgap}
\end{equation}
Thus, for an admissible subsolution $u,$
$\Theta_D(u) > 0$ is a violation of the mean value inequality.
In this notation, $D$ is a weak mean value set precisely
when $\Theta_D(u) = 0$ for every admissible $L$-harmonic
function $u,$ while $D$ is a strong mean value set precisely
when $\Theta_D(v) \le 0$ for every admissible $L$-subsolution
$v.$

We first relate the sign of the weak modified Schwarz
$L$-potential to the strong mean value property.
The argument is a divergence-form analogue of
\cite[Theorem 2.1(b)]{ShSj}; see also the duality
argument in \cite[Section 8]{GRooooos}.

\begin{lemma}[Mean-value gap duality and nonnegativity]
\label{nonneg}
Let $D \Subset \Omega$ be a weak mean value set for
$L$ at $x_0,$ where $\Omega$ is a bounded auxiliary
domain satisfying the standing auxiliary-domain assumptions.
For $f \in L^\infty(\Omega),$ let
$\psi_f \in W^{1,2}_0(\Omega)$ solve
Equation~\refeqn{psieq}\!\!. Then
\begin{equation}
\Theta_D(\psi_f)
  = -\int_\Omega W_D(y)f(y) \, dy.
\label{eq:dual}
\end{equation}
In particular, if $D$ is a strong mean value set, then
$W_D \ge 0$ almost everywhere in $\Omega.$
\end{lemma}

\begin{pf}
Green representation gives
$$\psi_f(x) = -\int_\Omega G(x,y)f(y) \, dy.$$
Using Fubini's theorem and the symmetry of $G,$ we obtain
\begin{align*}
\Theta_D(\psi_f)
 &= \psi_f(x_0)
    -\frac{1}{|D|}\int_D\psi_f(x)\,dx\\
 &= -\int_\Omega G(x_0,y)f(y)\,dy
    +\frac{1}{|D|}\int_D\int_\Omega
       G(x,y)f(y)\,dy\,dx\\
 &= -\int_\Omega
    \left(
      G(y,x_0)
      -\frac{1}{|D|}\int_D G(y,x)\,dx
    \right)f(y)\,dy\\
 &= -\int_\Omega W_D(y)f(y)\,dy,
\end{align*}
where the last equality follows from
Equation~\refeqn{potdef}\!\!.
This proves Equation~\refeqn{dual}\!\!.

If $D$ is strong and $f\ge0,$ then
$L\psi_f=f\ge0,$ so $\psi_f$ is an admissible
subsolution and hence $\Theta_D(\psi_f)\le0.$
Equation~\refeqn{dual} therefore gives
$$\int_\Omega W_D f\,dx\ge0$$
for every nonnegative $f\in L^\infty(\Omega),$
which implies $W_D\ge0$ almost everywhere.
\end{pf}

We next compare potentials associated with sets of the
same volume and the same point source.

\begin{theorem}[Comparison and uniqueness of potentials]
\label{SUR}
Let $D_1,D_2 \Subset \Omega$ be measurable sets with
$$\mathcal{V} := |D_1| = |D_2| \in (0,\infty).$$
Suppose that for a fixed
$x_0 \in \Omega, \ W_j \in W^{1,1}_{\mathrm{loc}}(\Omega)$ satisfy
$$LW_j = \frac{1}{\mathcal{V}}\bigrchi_{D_j} - \delta_{x_0}
  \quad\text{in }\Omega,
  \qquad
W_j = 0 \quad\text{a.e. in }\Omega\setminus D_j,$$
for $j = 1,2,$ and that
$U := W_1-W_2 \in W^{1,2}_0(\Omega).$
Then
\begin{equation}
-\mathcal{V} \int_{\Omega} A\nabla U \cdot \nabla U \, dx
  = \int_{D_1\setminus D_2}W_1 \, dx
    + \int_{D_2\setminus D_1}W_2 \, dx.
\label{eq:cmp}
\end{equation}
If $W_2 \ge 0$ a.e., then $W_1 \le W_2$ a.e.
If both potentials are nonnegative, then
$$W_1 = W_2,
  \qquad \bigrchi_{D_1} = \bigrchi_{D_2}
  \quad\text{a.e. in }\Omega.$$
\end{theorem}

\begin{pf}
Subtracting the distributional equations for
$W_1$ and $W_2$ gives
$$LU = \frac{1}{\mathcal{V}}
  \left(\bigrchi_{D_1}-\bigrchi_{D_2}\right).$$
The right-hand side belongs to $L^2(\Omega),$ so
the identity extends by density to tests in
$W^{1,2}_0(\Omega).$
Testing with $U$ and using the vanishing of each
potential outside its set gives Equation~\eqref{eq:cmp}.

Suppose that $W_2 \ge 0.$ Then
$U_+ := \max\{U,0\} \in W^{1,2}_0(\Omega)$
vanishes a.e. outside $D_1.$ Testing with $U_+$ gives
$$-\int_\Omega A\nabla U_+\cdot\nabla U_+ \, dx
  = \frac{1}{\mathcal{V}}\int_{D_1 \setminus D_2} U_+ \, dx
  \ge 0.$$
Uniform ellipticity gives $\nabla U_+ = 0$ a.e.
The zero extension of $U_+$ belongs to
$W^{1,2}(\R^n)$ and has zero weak gradient;
its square integrability therefore forces $U_+ = 0.$
This proves $W_1 \le W_2.$
If both potentials are nonnegative, interchanging
the indices gives $W_1 = W_2.$ The difference equation
then gives $\bigrchi_{D_1} = \bigrchi_{D_2}$ a.e.
\end{pf}

\begin{lemma}[The canonical comparison]
\label{cancomp}
Let $D$ be a weak mean value set for $L$ at $x_0,$ put
$\mathcal{V}:=|D|$ and $r:=\mathcal{V}^{1/n},$ and let
$$E:=D_r(x_0)$$
be the canonical set with the same center and volume.
Choose a bounded auxiliary domain $\Omega$ satisfying the
standing auxiliary-domain assumptions and compactly
containing $D$ and $E.$  Let $W_D$ be the weak modified
Schwarz $L$-potential of $D,$ let $W_E:=W_r$ be the
canonical height of $E,$ and define the \textit{canonical
comparison}
\begin{equation}
\begin{aligned}
\Upsilon_{D,E}
  := W_D-W_E
  &= (W_D-G)+u_r
     \in W^{1,2}_0(\Omega),\\
L\Upsilon_{D,E}
  &= \frac{1}{\mathcal{V}}
     \left(\bigrchi_D-\bigrchi_E\right),
\end{aligned}
\label{eq:can}
\end{equation}
where $G:=G(\cdot,x_0).$
The function $\Upsilon_{D,E}$ has a locally \Holder continuous
representative throughout $\Omega$ and satisfies
$$\Upsilon_{D,E}\le0 \qquad\text{in }\Omega,$$
as well as
$$\Upsilon_{D,E}=0
  \qquad\text{on }\Omega\setminus(D\cup E).$$
\end{lemma}

\begin{pf}
Theorem~\refthm{MVSYWMSP}and $W_E=G-u_r$ give the
first line of Equation~\refeqn{can}\!\!, while subtraction
of the potential equations gives the second.
Since $W_E\ge0,$ Theorem~\refthm{SUR}gives
$\Upsilon_{D,E}\le0$ almost everywhere.
The bounded right-hand side in Equation~\refeqn{can}gives
a locally \Holder continuous representative throughout
$\Omega,$ so the inequality holds everywhere.
Both potentials vanish outside their respective sets,
which gives
$\Upsilon_{D,E}=0$ on $\Omega\setminus(D\cup E).$
\end{pf}

\begin{corollary}[Extremal characterization of the canonical potential]
\label{maxpotential}
Fix a center $x_0$ and a volume $\mathcal{V} > 0,$ and set
$$r := \mathcal{V}^{1/n}, \qquad E := D_r(x_0).$$
Among all weak mean value sets $D$ for $L$ at $x_0$
with $|D| = \mathcal{V},$ the canonical height $W_E$ is the
unique maximum of the associated weak modified Schwarz
$L$-potentials in the pointwise order. More precisely,
on any common bounded auxiliary domain $\Omega$
as in Lemma~\refthm{cancomp}\!\!,
$$W_D = W_E \quad\text{on }\Omega\setminus\{x_0\}
  \qquad\Longleftrightarrow\qquad
  |D\triangle E| = 0.$$
If $|D\triangle E|>0,$ then
$$\Upsilon_{D,E} < 0
  \qquad \text{throughout }D\cup E.$$
Equivalently,
$$W_D < W_E
  \quad \text{on }(D\cup E)\setminus\{x_0\},
  \qquad
D\setminus E
  \subset \{x\in D:W_D(x)<0\}.$$
\end{corollary}

\begin{pf}
Lemma~\refthm{cancomp}gives $W_D \le W_E.$
If $W_D = W_E,$ then Equation~\refeqn{can}gives
$\bigrchi_D = \bigrchi_E$ almost everywhere.
Conversely, if $|D\triangle E|=0,$ the
Green-potential formula in Equation~\refeqn{potdef}gives
$W_D = W_E.$

Now suppose that $|D\triangle E| > 0.$
By Lemma~\refthm{cancomp}\!\!,
$$L\Upsilon_{D,E}
  = \frac{1}{\mathcal{V}}(1-\bigrchi_E)\ge0
  \qquad\text{in }D,$$
while $\Upsilon_{D,E}\le0.$
Since $D$ is connected, the strong maximum principle
\cite[Theorem 8.19]{GT} shows that either
$\Upsilon_{D,E}<0$ throughout $D$ or
$\Upsilon_{D,E}\equiv0$ in $D.$
In the latter case, the preceding equation gives
$|D\setminus E| = 0.$ Since $|D| = |E|,$ this would imply
$|D\triangle E| = 0,$ a contradiction.
Hence
$$\Upsilon_{D,E} < 0 \qquad\text{throughout }D.$$

If $x \in E\setminus D,$ then $x\ne x_0,$
$W_D(x) = 0,$ and $W_E(x) > 0,$ so
$\Upsilon_{D,E}(x) < 0.$
Thus $\Upsilon_{D,E}<0$ throughout $D\cup E.$
Finally, on $D \setminus E$ we have $W_E=0,$ so
$W_D = \Upsilon_{D,E}<0,$ proving the last assertion.
\end{pf}

The same energy calculation gives quantitative control
in terms of a mean-value gap. The normalized harmonic
mean value gap of \cite{CFLZ}, with harmonic tests replaced
by $L$-harmonic tests, vanishes on every weak mean value set.
Theorem~\refthm{WDNIS}shows that such a purely harmonic gap
cannot control the measure of the symmetric difference from
the canonical set in general. The quantity below instead
measures the largest mean-value gap among the normalized
family of subsolutions $\psi_f$ with $0 \le f \le 1.$

\begin{proposition}[Maximal mean-value gap and potential control]
\label{potdefect}
Let $D$ be a weak mean value set for $L$ at $x_0,$ set
$r := |D|^{1/n},$ and let $E := D_r(x_0).$
Choose a bounded auxiliary domain $\Omega$ satisfying
the standing auxiliary-domain assumptions and compactly
containing both sets. Let $W_D$ be the weak modified
Schwarz $L$-potential of $D,$ and let $W_E := W_r$ be
the canonical height of $E.$ For $f \in L^\infty(\Omega),$
let $\psi_f$ solve Equation~\refeqn{psieq}\!\!, and define
the maximal mean-value gap
\begin{equation}
\delta(D)
 := \sup_{\substack{f\in L^\infty(\Omega)\\0\le f\le 1}}
       \Theta_D(\psi_f)
  = \int_{\Omega} (W_D)_-\,dx, \qquad
  (W_D)_-:=\max\{-W_D,0\}.
\label{eq:maxgap}
\end{equation}
The second equality in Equation~\refeqn{maxgap}follows
from Lemma~\refthm{nonneg}\!\!, and the supremum is
attained by $f = \bigrchi_{\{W_D<0\}}.$
The value of $\delta(D)$ is independent of the choice of $\Omega.$

Let $\Upsilon_{D,E}$ be the canonical comparison from
Lemma~\refthm{cancomp}\!\!. Then
$$\lambda |D|\int_\Omega |\nabla\Upsilon_{D,E}|^2\,dx
  + \int_{E\setminus D}W_E\,dx
  \le \delta(D).$$
In particular,
$$\|\nabla\Upsilon_{D,E}\|_{L^2(\Omega)}
     \le
  \left(\frac{\delta(D)}{\lambda |D|}\right)^{1/2}.$$
\end{proposition}

\begin{pf}
Let $\Omega_1$ and $\Omega_2$ be two such auxiliary domains.
If admissible sources on $\Omega_1$ and $\Omega_2$ agree
almost everywhere on $D,$ their Dirichlet solutions differ
by an $L$-harmonic function on $D.$ This difference is
integrable on $D$ since $D \Subset \Omega_1 \cap \Omega_2.$
The weak mean value property therefore gives the same
mean-value gap for both solutions. Every $f \in L^\infty(D)$
with $0 \le f \le 1$ extends by zero to either domain,
so the supremum in Equation\refeqn{maxgap}is the same for both domains.

Applying Equation~\refeqn{cmp}to the canonical comparison
$\Upsilon_{D,E}$ from Lemma~\refthm{cancomp}and using
uniform ellipticity gives
\begin{equation}
\begin{aligned}
\lambda |D|\int_\Omega |\nabla\Upsilon_{D,E}|^2\,dx
+\int_{E\setminus D}W_E\,dx
&\le
|D|\int_\Omega
  A\nabla\Upsilon_{D,E}\cdot\nabla\Upsilon_{D,E}\,dx
+\int_{E\setminus D}W_E\,dx\\
&=
-\int_{D\setminus E}W_D\,dx\\
&\le
\int_{D\setminus E}(W_D)_-\,dx
\le \delta(D).
\end{aligned}
\label{eq:TwoFer}
\end{equation}
Since $W_E\ge0,$ dropping its integral gives
the gradient bound.
\end{pf}

\begin{proposition}[Symmetric difference control]
\label{setdefect}
Under the hypotheses and notation of
Proposition~\ref{potdefect}, define
$m_E(t):=
  \bigl|\{x \in E : 0 < W_E(x)\le t\}\bigr|,$ for $t > 0.$
Then
$$|D \triangle E| \le
  2m_E(t) + \frac{2\delta(D)}{t} \qquad\text{for every }t > 0.$$
\end{proposition}

\begin{pf}
For $t > 0,$ split $E\setminus D$ according to whether
$W_E \le t$ or $W_E > t.$ Since $W_E > 0$ almost everywhere
in $E,$ Equation~\refeqn{TwoFer}gives
$$|E\setminus D| \le m_E(t)
  + \frac{1}{t} \int_{E \setminus D} W_E \, dx
  \le m_E(t) + \frac{\delta(D)}{t}.$$
Finally, $|D| = |E|$ gives
$|D \triangle E| = 2|E \setminus D|.$
\end{pf}

\section{Uniqueness and Classification}
\label{TUT}

\subsection{Classification of strong mean value sets}

For the Laplacian, the characterization by nonnegativity
of the modified Schwarz potential is classical; see
\cite[Proposition 8.3]{S} and \cite[Section 4.4]{Sh}.
We now prove Theorem~\refthm{openreps}\!\!. Throughout the proof,
$E := D_r(x_0)$ is the canonical set with the same center
and volume as $D.$

\begin{pf}
For $(1)\Rightarrow(2),$ a strong mean value set is
weak, and Lemma~\refthm{nonneg}gives nonnegativity
of its potential.

Suppose that $(2)$ holds. Choose a common bounded auxiliary
domain $\Omega$ as in Lemma~\refthm{cancomp}and apply that lemma
to obtain
$$\Upsilon_{D,E} = W_D - W_E \in W^{1,2}_0(\Omega).$$
Since both $W_D$ and $W_E$ are nonnegative,
Theorem~\refthm{SUR}gives
$$W_D = W_E \quad\text{a.e. in }\Omega,
  \qquad |D\triangle E| = 0.$$
Both potentials are continuous away from $x_0,$
so their equality holds pointwise there.
Since $W_D$ vanishes at every point outside $D,$
whereas $W_E > 0$ on $E\setminus\{x_0\},$ we obtain
$E \subset D;$ here we also use $x_0 \in D.$
Moreover, $D\setminus\overline{E}$ is an open set
of measure zero and is therefore empty.
Since $D$ is open, this proves $(3).$

Conversely, suppose that $(3)$ holds, so that
$E \subset D \subset \operatorname{int}\overline{E}.$
Then $D\setminus E \subset \partial E,$ and
Theorem~\refthm{KeyMVT}gives $|\partial E| = 0.$
Thus $D$ and $E$ agree almost everywhere.
Let $v \in L^1(D)\cap W^{1,2}_{\mathrm{loc}}(D)$ satisfy
$Lv \ge 0$ weakly in $D.$
Its restriction to $E$ is an integrable subsolution
there, so Theorem~\refthm{KeyMVT}\!\!(3) gives
$$v(x_0)
  \le \frac{1}{|E|}\int_E v(x) \, dx
  = \frac{1}{|D|}\int_D v(x) \, dx.$$
Hence $D$ is strong.

For a weak mean value set $D,$ Equation\refeqn{maxgap}gives
$\delta(D) = 0$ if and only if $W_D \ge 0$ almost everywhere.
By continuity away from $x_0,$ this is equivalent to pointwise
nonnegativity there. Thus $(2)$ and $(4)$ are equivalent.

Finally, $E$ and $\operatorname{int}\overline{E}$
are bounded open sets containing $x_0,$ and both
have volume $r^n$ because their difference is
contained in $\partial E.$
Applying the characterization to these two sets
proves the remaining assertions.
\end{pf}

\begin{corollary}[Almost-everywhere uniqueness]
\label{AMVSAaNSfOP}
Let $n \ge 2,$ and let $L := \div(A\nabla\cdot)$ on $\R^n,$
where $A$ is bounded, measurable, symmetric, and uniformly
elliptic. If $D$ is a strong mean value set for $L$ at
$x_0 \in D,$ then
\begin{equation}
  \bigl|D \triangle D_{|D|^{1/n}}(x_0)\bigr| = 0.
\label{eq:main}
\end{equation}
\end{corollary}

\begin{pf}
Take $r := |D|^{1/n}$ in Theorem~\refthm{openreps}\!\!.
\end{pf}

\subsection{Point deletion}

\begin{proposition}[Mean-value gap and point deletion]
\label{pointdefect}
Let $S$ be a bounded open strong mean value set for $L$
at $x_0,$ and choose $r > 0$ so that $r^n = |S|.$
Set $E := D_r(x_0),$
and let $W_E$ be the associated canonical height function.
Choose a common sufficiently large bounded smooth ambient
domain $\Omega$ compactly containing $S$ and $E,$ and let
$G$ be its positive Dirichlet Green function, normalized by
$L G(\cdot,y) = -\delta_y.$

For any $y \in S \setminus \{x_0\},$ set
$S_y := S \setminus \{y\}$ and $g_y(x) := G(x,y).$
Then $g_y$ is an admissible integrable $L$-harmonic function
on $S_y,$ and
$$\Theta_{S_y}(g_y) = W_E(y).$$
Moreover, the following conditions are equivalent:
\begin{enumerate}
\item $S_y$ is a weak mean value set for $L$ at $x_0.$
\item $S_y$ is a strong mean value set for $L$ at $x_0.$
\item $W_E(y) = 0.$
\item $y \in S \setminus E.$
\end{enumerate}
In particular, deleting any point of
$E \setminus \{x_0\}$ from $S$ destroys even the weak
mean value property.
\end{proposition}

\begin{pf}
The Green function $g_y$ is locally integrable across its
pole and belongs to $W^{1,2}_{\mathrm{loc}}(S_y).$
Since its pole is outside $S_y,$ it satisfies
$L g_y = 0$ weakly in $S_y.$
Thus it is an admissible test function on $S_y.$

Let $W_S$ be the potential defined by
Equation~\refeqn{potdef}with $D = S.$
Theorem~\refthm{openreps}gives
$$W_S = W_E
  \qquad\text{on }\Omega\setminus\{x_0\}.$$
Since $|S_y| = |S|,$ deleting $y$ does not change the
integral of $g_y.$
Using the symmetry of the Green function, we obtain
\begin{equation}\label{eq:GapIdent}
  \begin{aligned}
  \Theta_{S_y}(g_y)
   &= g_y(x_0)
    - \frac{1}{|S_y|}\int_{S_y} g_y(x) \, dx\\
   &= G(x_0,y)
    - \frac{1}{|S|}\int_S G(x,y) \, dx \\
   &= G(y,x_0)
    - \frac{1}{|S|}\int_S G(y,x) \, dx \\
   &= W_S(y) \\
   &= W_E(y).
  \end{aligned}
\end{equation}

If $S_y$ is a weak mean value set, then
$\Theta_{S_y}(g_y) = 0,$ since $g_y$ is an admissible
$L$-harmonic test function. Equation~\eqref{eq:GapIdent}
therefore gives $W_E(y) = 0.$
Since $y \ne x_0$ and $E$ is the canonical positivity set,
this is equivalent to $y \in S \setminus E.$

Conversely, suppose that $y \in S\setminus E.$
Theorem~\refthm{openreps}gives
$$E \subset S_y \subset S
  \subset \operatorname{int}\overline{E}.$$
Since $S_y$ is open, contains $x_0,$ and has volume
$r^n,$ another application of that theorem shows
that $S_y$ is strong, hence also weak.

Finally, if $y \in E \setminus \{x_0\},$ then
Equation~\eqref{eq:GapIdent} gives
$$\Theta_{S_y}(g_y) = W_E(y)>0.$$
Thus $g_y$ explicitly gives a mean-value violation
and witnesses the failure of the weak mean value property.
\end{pf}

\begin{remark}[Deleting a point]
\label{deletepoint}
Deleting a set of measure zero leaves integrals unchanged
but can enlarge the class of admissible test functions.
For example, in $\R^2,$ let
$$a = (1/2,0), \ \ \
  D = B_1(0)\setminus\{a\}, \ \ \
  u(x) = -\log|x-a|.$$
Then $u$ is smooth and harmonic in $D$ and belongs to $L^1(D).$
For $0 < r < 1$ with $r \ne |a|,$ its average over
$\partial B_r(0)$ is $-\log\max\{r,|a|\}.$ Integrating in polar
coordinates gives
$$\frac{1}{|D|}\int_D u(x) \, dx
  = -2\int_0^1 r\log\max\{r,|a|\} \, dr
  = \frac{1 - |a|^2}{2} = \frac{3}{8} < \log 2 = u(0),$$
so $\Theta_D(u)>0.$ Thus $D$ is not even a weak mean
value set for the Laplacian at $0,$ despite agreeing
almost everywhere with $B_1(0).$ The point is that our
definition permits test functions defined on $D$ itself;
they need not extend harmonically across the deleted point.
In particular, the normalized harmonic mean value gap
of \cite{CFLZ} is positive for this punctured ball,
although its symmetric difference from $B_1(0)$ has measure
zero.
\end{remark}

\subsection{Conditions under which weak implies strong}

The next result is a divergence-form analogue of the
connected-exterior criterion for harmonic and subharmonic
balls in \cite[Corollary 3.4]{ShSj}.
In the planar weighted setting, Hedenmalm and Shimorin
\cite[Proposition~2.5]{HS02} prove that if the Hele--Shaw
domain is a Jordan domain, the harmonic mean value identity
determines it up to a set of area zero.

\begin{proposition}[A connected-exterior criterion]
\label{extstrong}
Fix a center $x_0$ and a canonical parameter $r > 0,$
set $E := D_r(x_0),$ and suppose that
$\R^n \setminus \overline{E}$ is connected.
Then every weak mean value set $D$ for $L$ at $x_0$
with $|D| = r^n$ is a strong mean value set.
Consequently, the open weak and strong mean value sets
with this center and volume are precisely the open
sets satisfying
$$E \subset D \subset \operatorname{int}\overline{E}.$$
Equivalently, if there exists a weak mean value set for $L$
at $x_0$ with volume $r^n$ which is not strong, then
$\R^n \setminus \overline{D_r(x_0)}$ must be disconnected.
\end{proposition}

\begin{pf}
Let $D$ be such a weak mean value set, put
$\mathcal{V} := |D| = |E|,$ and choose a sufficiently large
ball $\Omega$ compactly containing both sets.
Write $W_E := W_r,$ and let
$\Upsilon_{D,E}$ be the canonical comparison from
Lemma~\refthm{cancomp}\!\!.
Then
$$\Upsilon_{D,E} \le 0
  \qquad \text{throughout }\Omega.$$

Set
$$O := \Omega \setminus \overline{E}.$$
Because $\Omega$ is a ball compactly containing $E,$
the connectedness of $\R^n \setminus \overline{E}$
implies that $O$ is connected.
Indeed, every point of $O$ can be joined within $O$
to an outer annulus in $\Omega$ disjoint from
$\overline{E},$ and that annulus is connected
because $n \ge 2.$

On $O,$ the canonical height vanishes, so
$$W_D = \Upsilon_{D,E} \le 0,
  \ \ \ LW_D = \frac{1}{\mathcal{V}}\bigrchi_D \ge 0.$$
Moreover, $W_D = 0$ in a nonempty open collar of
$\partial\Omega,$ since $D$ is compactly contained
in $\Omega.$
The strong maximum principle for weak subsolutions
therefore gives
$$W_D \equiv 0 \qquad \text{in } O;$$
see \cite[Theorem 8.19]{GT}.
Its equation now implies $|D \cap O| = 0.$
Since $D \cap O$ is open, it is empty.
Consequently,
$$D \subset \overline{E}.$$

Theorem~\ref{KeyMVT} gives $|\partial E| = 0,$ so
$|D \setminus E| = 0.$
Together with $|D| = |E|,$ this yields
$$|D \triangle E| = 0.$$
The Green representations on the common ambient domain
then give
$$W_D = W_E \ge 0
  \ \ \ \text{on } \Omega \setminus \{x_0\}.$$
Theorem~\refthm{openreps}now gives both the strong
mean value property and the asserted classification.
\end{pf}

\begin{corollary}[Enclosure under the standing coefficient assumptions]
\label{roughenc}
Let $D$ be a weak mean value set for $L$ at $x_0$ which is
not strong, and set $E := D_{|D|^{1/n}}(x_0).$
Then $\partial D \subset \overline{E},$ and every connected
component of $\R^n \setminus \overline{E}$ that meets $D$
is contained in $D.$ At least one such component exists.
\end{corollary}

\begin{pf}
Choose a ball $\Omega$ compactly containing $D$ and $E.$
If $z \in \partial D \setminus \overline{E},$ take a ball
$B$ centered at $z$ with $B \Subset \Omega \setminus \overline{E}.$
Lemma~\refthm{cancomp}gives $W_D \le W_E = 0$ in $B,$ and
$W_D(z) = 0$ since $z \notin D.$ As
$LW_D = |D|^{-1}\bigrchi_D \ge 0$ in $B,$ the strong maximum
principle gives $W_D \equiv 0$ in $B.$ This contradicts the
equation on the nonempty open set $B \cap D.$ Thus
$\partial D \subset \overline{E}.$

Every connected component of $\R^n \setminus \overline{E}$
therefore avoids $\partial D,$ so any such component meeting
$D$ lies in $D.$ Finally, the last paragraph of the proof of
Proposition~\refthm{extstrong}shows that $D \subset \overline{E}$
would imply that $D$ is strong. Hence $D \setminus \overline{E}$
is nonempty, and the component containing any of its points
is bounded since it lies in $D.$
\end{pf}

\begin{corollary}[Two classes for which weak implies strong]
\label{goodclasses}
In the two cases below, every weak mean value set of
volume $r^n$ at the indicated center coincides exactly
with the displayed canonical set and is therefore strong.
\begin{enumerate}
\item
Let $B$ be a constant symmetric positive definite
matrix and let $\Psi:\R^n \longrightarrow \R^n$
be a smooth bi-Lipschitz diffeomorphism with
$\det D\Psi \equiv J > 0.$ Define
$$A(\Psi(y))
  := \frac{D\Psi(y)B D\Psi(y)^T}{J}.$$
Then, at every center $x_0 = \Psi(y_0),$
$$D_r(x_0)
  = \Psi\left(y_0 + B^{1/2}B_{R_{r,J}}(0)\right),
  \qquad
  R_{r,J} := \frac{r}
  {(J\omega_n\sqrt{\det B})^{1/n}}.$$
Here $B^{1/2}$ denotes the positive definite square root
of $B.$  In particular, taking $\Psi = \operatorname{id}$
and $J = 1$ gives the constant-matrix case
$A \equiv B,$ for which the canonical sets are ellipsoids.

\item
If $A(x) = a(|x-z_0|)I,$ where $a$ is measurable and
$0 < \lambda \le a(t) \le \Lambda < \infty$
for almost every $t > 0,$ then
$$D_r(z_0) = B_{r/\omega_n^{1/n}}(z_0).$$
\end{enumerate}
\end{corollary}

\begin{pf}
For (1), the bi-Lipschitz hypothesis ensures uniform
ellipticity of $A.$
Let $D$ be a weak mean value set of volume $r^n$
at $x_0,$ and set
$$F(z) := \Psi(y_0 + B^{1/2}z).$$
Under the change of variables $x=F(z),$ the coefficient
matrix $A$ becomes $\sqrt{\det B}\,I$ in the $z$-coordinates.
Since multiplication by a positive constant does not change
the homogeneous equation, admissible $L$-harmonic functions
on $D$ correspond to integrable harmonic functions on
$F^{-1}(D);$ see \cite[Section 2.2]{KSVW}.
Since $\det DF \equiv J\sqrt{\det B}$ is constant,
normalized Lebesgue averages are preserved.
Thus $F^{-1}(D)$ is a weak mean value set for the
Laplacian at $0,$ of volume
$r^n/(J\sqrt{\det B}).$
Theorem~\ref{KurThm} identifies it with
$B_{R_{r,J}}(0).$
Applying this conclusion also to the canonical set
gives the asserted exact equality and hence the
strong mean value property.

For (2), rotational invariance and uniqueness of the
Green function and obstacle solution in a sufficiently
large ball centered at $z_0$ make $E := D_r(z_0)$
invariant under rotations about $z_0.$
Since $E$ is bounded, open, connected
(Remark~\ref{Ambient}), and contains $z_0,$
it is a ball; its volume $r^n$ determines its radius.
Its open exterior is connected, so
Proposition~\ref{extstrong} gives
$$E \subset D \subset \operatorname{int}\overline{E} = E$$
for every weak mean value set $D$ with this center
and volume. Thus $D = E$ and is strong.
\end{pf}
\section{Weak Mean Value Sets That Are Not Strong}
\label{sharpness}

\subsection{A smooth counterexample in every dimension}

Theorem~\refthm{WDNIS}shows that the strong mean value hypothesis
in Corollary~\refthm{AMVSAaNSfOP}is essential. A weak mean value
set can differ on a set of positive measure from the canonical
set with the same center and volume. This contradicts the
uniqueness assertion in our earlier preprint \cite{ABthepain}.

We prove Theorem~\refthm{WDNIS}after establishing two radial
ingredients. The first shows that the criterion in the second
is nonvacuous in every dimension.

\begin{lemma}[Radial maps satisfying an integral criterion]
\label{radexist}
For every $n \ge 2,$ there exists an orientation-preserving smooth
radial diffeomorphism
$\Psi : \R^n \longrightarrow \R^n,$ equal to the identity outside a
compact subset of $B_1(0),$ such that, for
\begin{equation}
  A(\Psi(y)) :=
  \frac{D\Psi(y)D\Psi(y)^T}{\det D\Psi(y)},
\label{eq:Aconduct}
\end{equation}
the radial eigenvalue $\alpha(r)$ of $A(x)$ at $|x| = r$ satisfies
$$\int_0^1 \frac{r}{\alpha(r)} \, dr > \frac{n + 1}{2}.$$
\end{lemma}

\begin{pf}
Fix $0 < a < \frac{1}{8},$ and define the strictly increasing
piecewise-linear function
$$\overline{s}_a(t) :=
  \begin{cases}
  \dfrac{t}{2a},
      &0 \le t \le a,\\[9pt]
  \dfrac{1}{2} + \dfrac{t - a}{3 - 4a},
      &a \le t \le \dfrac{3}{4},\\[7pt]
  t,
      &t \ge \dfrac{3}{4}.
  \end{cases}$$
Smooth $\overline{s}_a$ in the intervals
$(a - a/10,a + a/10)$ and
$(3/4 - a/10,3/4 + a/10),$ while making sure to have a
strictly positive derivative there,
and denote the resulting smooth function (defined as
$\overline{s}_a$ outside these intervals) by $s_a.$ In particular,
$s_a(t) \ge 1/2$ for $2a \le t \le 1/2.$
Define
$$\Psi_a(y) :=
  \begin{cases}
  \dfrac{s_a(|y|)}{|y|}y, & y \ne 0,\\[6pt]
  0, & y = 0.
  \end{cases}$$
Because $s_a$ is linear near $0,$ strictly increasing, and equal to
the identity near infinity, $\Psi_a$ is an orientation-preserving
smooth radial diffeomorphism equal to the identity outside a compact
subset of $B_1(0).$

Let $A_a$ be the matrix field defined by
Equation\refeqn{Aconduct}with $\Psi = \Psi_a,$ and let $\alpha_a$
be the radial eigenvalue of $A_a.$ Fix $y = te,$ where $|e| = 1.$
Then the radial vector $e$ and every unit tangential vector
$\tau \perp e$ are eigenvectors of $D\Psi_a(te),$ with
$$D\Psi_a(te)e = s_a'(t)e,
  \qquad
  D\Psi_a(te)\tau = \frac{s_a(t)}{t}\tau.$$
Thus $D\Psi_a(te)$ is diagonal in an orthonormal radial-tangential
basis, and hence
$$\det D\Psi_a(te)
  = s_a'(t)\left(\frac{s_a(t)}{t}\right)^{n - 1}.$$
Applying Equation\refeqn{Aconduct}to the radial vector $e$ gives
$$A_a(s_a(t)e)e
  = \frac{s_a'(t)^2}
          {s_a'(t)\left(s_a(t)/t\right)^{n - 1}}e
  = s_a'(t)
     \left(\frac{t}{s_a(t)}\right)^{n - 1}e,$$
and thus,
$$\alpha_a(s_a(t))
  = s_a'(t)\left(\frac{t}{s_a(t)}\right)^{n - 1}.$$
Changing variables $r = s_a(t)$ therefore gives
$$\int_0^1 \frac{r}{\alpha_a(r)} \, dr
  = \int_0^1 \frac{s_a(t)}{s_a'(t)}
       \left(\frac{t}{s_a(t)}\right)^{1 - n} s_a'(t) \, dt
  = \int_0^1 s_a(t)^n t^{1 - n} \, dt
  \ge \frac{1}{2^n}\int_{2a}^{1/2} t^{1 - n} \, dt.$$
The last integral tends to infinity as $a \downarrow 0$ for every
$n \ge 2.$ Choosing $a$ sufficiently small proves the lemma.
\end{pf}

\begin{proposition}[A radial criterion in every dimension]
\label{radcrit}
Fix $n \ge 2,$ $p \in \R^n,$ and $R > 0.$ Let
$\Psi : \R^n \longrightarrow \R^n$ be an orientation-preserving
smooth radial diffeomorphism about $p,$ equal to the identity
outside a compact subset of $B_R(p).$
Write $\Psi(p + te) = p + s(t)e$ for $t > 0,\ |e| = 1,$
and define $A$ by the standard conductivity transformation
$$A(\Psi(y)) :=
  \frac{D\Psi(y)D\Psi(y)^T}{\det D\Psi(y)},
  \qquad L := \div(A\nabla\cdot).$$
Let $\alpha(r)$ denote the radial eigenvalue of $A(x)$ at points
$x$ satisfying $|x - p| = r,$ so that
$\alpha(s(t)) = s'(t)\left(\frac{t}{s(t)}\right)^{n - 1}.$
Assume
\begin{equation}
\label{eq:radcriterion}
\int_0^R \frac{r}{\alpha(r)} \, dr
  > \frac{n + 1}{2}R^2,
  \qquad \text{and set } \quad
  R_{\mathrm{crit}}
  := \left(
       \frac{2}{n}\int_0^R \frac{r}{\alpha(r)} \, dr
       - \frac{R^2}{n}
     \right)^{1/2}.
\end{equation}
Then $R_{\mathrm{crit}} > R,$ and for every
$R < K < R_{\mathrm{crit}},$ there exists $d_0 > K$ such that,
for every $x_0$ with $d := |x_0 - p| > d_0$ and every
$\varrho \in [d + R, d + K],$
the ball $D_\varrho := B_\varrho(x_0) \supset B_R(p)$ is a
weak mean value set for $L$ at $x_0$ but is not strong.
\end{proposition}

\begin{pf}
The hypotheses on $\Psi$ ensure that $A$ is smooth, symmetric, and
uniformly elliptic. Put $\rho := \det D\Psi.$ The radial change of
variables gives
$\rho(p + te)
  = s'(t)\left(\frac{s(t)}{t}\right)^{n - 1}.$

The Jacobian identity gives the mass-transport relation
$$
\Psi_{\#}\bigl(\rho(y)\,dy\bigr) = dx;
\quad\text{equivalently,}
\quad \forall \text{ measurable }E \subset \R^n, \quad
|E| = \int_{\Psi^{-1}(E)} \rho(y) \, dy.
$$
Thus $\Psi$ is the monotone radial transport, familiar from optimal
transport, which sends the weighted measure $\rho(y) \, dy$ to
Lebesgue measure. The standard conductivity transformation gives,
for smooth $v,$
\begin{equation}
\label{eq:condtrans}
Lv(\Psi(y))
  = \frac{\Delta(v \circ \Psi)(y)}{\rho(y)},
\end{equation}
and the same identity holds for weak equations; see
\cite[Section 2.2]{KSVW} and, in the plane,
\cite[Section 1]{ALP}.

Let $\phi$ be the solution of
\begin{equation}
\label{eq:phiBVP}
\begin{cases}
\Delta \phi = \rho - 1 & \text{in } B_R(p),\\
\phi = 0 & \text{on } \partial B_R(p).
\end{cases}
\end{equation}
By rotational invariance and uniqueness, $\phi$ is radial about
$p.$ Write $\phi(p + te) = \Phi(t)$ for $|e| = 1.$
For $0 < t < R,$ we have $\Psi(B_t(p)) = B_{s(t)}(p).$
The divergence theorem and the change of variables $x = \Psi(y)$ give
$$\begin{aligned}
n\omega_n t^{n - 1} \Phi'(t)
  &= \int_{B_t(p)} (\rho - 1) \, dy\\
  &= |B_{s(t)}(p)| - |B_t(p)|
   = \omega_n \bigl( s(t)^n - t^n \bigr).
\end{aligned}$$
Since $s(t) = t$ near $R,$ we have $\Phi'(t) = 0$ there.
Together with $\Phi(R) = 0,$ this shows that $\phi$ vanishes
in an outer annulus of $B_R(p).$ Its zero extension therefore
satisfies
$$\phi \in C_c^\infty(B_R(p)), \qquad
  \Delta \phi = \rho - 1 \quad \text{in } \R^n.$$
Solving the flux identity for $\Phi'$ and integrating from
$t$ to $R$ gives
$$\phi(p + te)
  = -\frac{1}{n}\int_t^R
       \frac{s(\tau)^n - \tau^n}{\tau^{n - 1}} \, d\tau
  \qquad (0 \le t < R).$$
In particular, using $r = s(t)$ and the formula for $\alpha,$
\begin{equation}
\label{eq:phicrit}
\begin{aligned}
\phi(p)
  &= \frac{R^2}{2n}
     - \frac{1}{n}\int_0^R s(t)^n t^{1 - n} \, dt\\
  &= \frac{R^2}{2n}
     - \frac{1}{n}\int_0^R \frac{r}{\alpha(r)} \, dr\\
  &= -\frac{R_{\mathrm{crit}}^2}{2}
   < -\frac{R^2}{2}.
\end{aligned}
\end{equation}

Fix $K$ with $R < K < R_{\mathrm{crit}}.$ For $x_0$ with
$d := |x_0 - p| > K$ and $\varrho \in [d + R, d + K],$ put
$k := \varrho - d,$ and $D := B_\varrho(x_0) = B_{d + k}(x_0).$
Then $R \le k \le K,$ and $B_R(p) \subset D.$ Consequently,
$\Psi$ preserves $D$ and is the identity near $x_0.$

We first verify the weak mean value property. Let
$v \in L^1(D) \cap W^{1,2}_{\mathrm{loc}}(D)$ satisfy $Lv = 0$ weakly in
$D,$ and put $h := v \circ \Psi.$ Since $\rho$ is bounded above
and below by positive constants, $h \in L^1(D).$
Equation\refeqn{condtrans}shows that $h$ is harmonic in $D.$
Since $\phi \in C_c^\infty(D),$ integration by parts gives
$$\begin{aligned}
\int_D v(x) \, dx
  &= \int_D h(y)\rho(y) \, dy\\
  &= \int_D h \, dy + \int_D h\Delta\phi \, dy\\
  &= |D|h(x_0) = |D|v(x_0).
\end{aligned}$$
Thus $D$ is a weak mean value set for $L$ at $x_0.$

We next show that $D$ is not strong. Let $P_D$ be the
unnormalized canonical height for the Laplacian on the ball $D,$
characterized distributionally by
\begin{equation}
\label{eq:PDBVP}
\Delta P_D = \bigrchi_D - |D|\delta_{x_0}
  \quad \text{in } \R^n,
  \qquad P_D = 0 \quad \text{outside } D.
\end{equation}
This function divided by $|D|$ is used in Caffarelli's obstacle
problem proof of the Mean Value Theorem in the third section of
the Fermi Lectures \cite{C2}.  By rotational invariance and
uniqueness, $P_D$ is radial about $x_0.$
Writing $P_D(x) = \mathcal{P}(|x - x_0|),$ integration of
Equation\refeqn{PDBVP}over $B_r(x_0)$ gives, for
$0 < r < d + k,$
$$n\omega_n r^{n - 1}\mathcal{P}'(r)
  = \omega_n r^n - \omega_n(d + k)^n.$$
Since $\mathcal{P}(d + k) = 0,$ we obtain
$$\mathcal{P}(r)
  = \frac{1}{n}\int_r^{d + k}
      \left((d + k)^n t^{1 - n} - t\right) \, dt > 0
  \qquad (0 < r < d + k).$$
In particular,
\begin{equation}
\label{eq:PDasymp}
\begin{aligned}
P_D(p)
  &= \frac{1}{n}\int_d^{d + k}
       \left((d + k)^n t^{1 - n} - t\right) \, dt\\
  &= \frac{1}{n}\int_0^k
       \left(\frac{(d + k)^n}{(d + t)^{n - 1}} - (d + t)\right)
       \, dt
   \longrightarrow \frac{k^2}{2}
     \quad \text{as } d \longrightarrow \infty.
\end{aligned}
\end{equation}
Factoring the difference of $n$th powers, we can write the last
integrand, including the factor $1/n,$ as
$$\frac{k - t}{n}\sum_{j = 0}^{n - 1}
    \left(\frac{d + k}{d + t}\right)^j.$$
For $0 \le t \le k \le K,$ the ratio lies between $1$ and
$1 + K/d.$ Thus the integrand converges uniformly to $k - t,$
and integration from $0$ to $k$ gives the stated limit uniformly
for $R \le k \le K.$

Set $Q_D := P_D + \phi.$ Since
$\operatorname{supp}\phi \subset B_R(p) \subset D,$ we have
\begin{equation}
\label{eq:QDeqn}
\Delta Q_D = \rho\bigrchi_D - |D|\delta_{x_0}
  \quad \text{in } \R^n,
  \qquad Q_D = 0 \quad \text{outside } D.
\end{equation}
Define
$$V_D(x) := \frac{Q_D(\Psi^{-1}(x))}{|D|}.$$
The weak form of Equation\refeqn{condtrans}gives
\begin{equation}
\label{eq:VDeqn}
LV_D = \frac{1}{|D|}\bigrchi_D - \delta_{x_0},
  \qquad V_D = 0 \quad \text{outside } D.
\end{equation}
On a sufficiently large bounded smooth ambient domain $\Omega,$
let $G(\cdot,x_0)$ be the positive Dirichlet Green function for
$L.$ Near $x_0,$ the map $\Psi$ is the identity, $A = I,$ and
$\phi = 0,$ so both $V_D$ and $G(\cdot,x_0)$ differ from the same
fundamental-solution singularity by smooth functions. Hence
$V_D - G(\cdot,x_0) \in W^{1,2}_0(\Omega).$
Equation\refeqn{VDeqn}and Theorem~\refthm{MVSYWMSP}therefore
identify $V_D$ with the weak modified Schwarz $L$-potential of $D.$

Finally, Equation\refeqn{phicrit}and $k \le K$ give
$$\frac{k^2}{2} + \phi(p)
  \le \frac{K^2 - R_{\mathrm{crit}}^2}{2} < 0.$$
By Equation\refeqn{PDasymp}and the uniformity just established,
there exists $d_0 > K,$
independent of the direction of $x_0 - p$ and of
$k \in [R, K],$ such that
$P_D(p) + \phi(p) < 0$ whenever $d > d_0.$
Since $\Psi(p) = p,$ we have $V_D(p) < 0.$
Theorem~\refthm{WisSiffPisP}shows that $D$ is not strong.
\end{pf}

We can now prove Theorem~\refthm{WDNIS}\!\!.

\begin{pf}
Choose $\Psi$ from Lemma~\refthm{radexist}and define $A$ and $L$
as in Proposition~\refthm{radcrit}with $p = 0$ and $R = 1.$
The proposition gives a ball
$D := B_\varrho(x_0)$
which is a weak mean value set for $L$ at $x_0$ but is not strong.
Let
$T(z) := x_0 + \varrho z,$
and define
$$\widehat{\Psi} := T^{-1} \circ \Psi \circ T,
  \qquad \widehat{A}(z) := A(T(z)),
  \qquad \widehat{L} := \div(\widehat{A}\nabla\cdot).$$
Then
$$\widehat{A}(\widehat{\Psi}(z))
  = \frac{D\widehat{\Psi}(z)D\widehat{\Psi}(z)^T}
         {\det D\widehat{\Psi}(z)}.$$
The affine change of variables $x = T(z)$ preserves normalized
averages and the signs of the corresponding differential
inequalities. Hence $B_1(0)$ is a weak mean value set for
$\widehat{L}$ at $0$ but is not strong.

Moreover, $\widehat{\Psi}$ is the identity near $0$ and outside a
compact subset of $B_1(0).$ Thus $\widehat{A} = I$ outside a compact
subset of $B_1(0).$ Finally,
$\det \widehat{A}(\widehat{\Psi}(z))
  = \bigl(\det D\widehat{\Psi}(z)\bigr)^{2 - n},$
so $\det \widehat{A} = 1$ when $n = 2.$ Renaming
$\widehat{\Psi}$ and $\widehat{A}$ as $\Psi$ and $A$ proves the
theorem.
\end{pf}

\begin{remark}[A smooth witness for failure of the strong property]
\label{smoothwitness}
The proof of Proposition~\refthm{radcrit}also gives a globally
smooth strict subsolution witnessing failure of the strong property.
Indeed, for one of the balls $D$ produced there, choose a nonzero
$f \in C_c^\infty(\{Q_D < 0\})$ with $f \ge 0.$ Let $\Gamma_n$ be
the fundamental solution normalized by $\Delta\Gamma_n = \delta_0,$
and, for $\varepsilon > 0,$ set
$$u(y) := (\Gamma_n * f)(y) + \varepsilon |y|^2.$$
Then $u \in C^\infty(\R^n)$ and
$$\Delta u = f + 2n\varepsilon > 0.$$
Because $Q_D$ is compactly supported and
Equation\refeqn{QDeqn}holds\!\!, we have
$$\begin{aligned}
\int_D u\rho \, dy - |D|u(x_0)
  &= \int_{\R^n}Q_D\Delta u \, dy\\
  &= \int_{\R^n}Q_D f \, dy
     + 2n\varepsilon\int_{\R^n}Q_D \, dy < 0
\end{aligned}$$
for sufficiently small $\varepsilon > 0.$ Set
$v := u \circ \Psi^{-1}.$ Equation\refeqn{condtrans}gives
$Lv > 0$ in $\R^n,$ while the change of variables gives
$$\frac{1}{|D|}\int_D v(x) \, dx - v(x_0) < 0.$$
Thus the failure is witnessed by a globally smooth strict
subsolution. The affine rescaling in the proof of
Theorem~\refthm{WDNIS}preserves this property for the unit-ball
counterexample.
\end{remark}

\subsection{Sufficient conditions on radial coefficients}

In the plane, Corollary~\refthm{plancrit}expresses the radial
criterion in Proposition~\refthm{radcrit}directly in terms of the
coefficients and gives explicit ellipticity bounds.

\begin{corollary}[A planar coefficient condition]
\label{plancrit}
Fix $p \in \R^2$ and $R > 0.$ Let $a$ be a positive
smooth function on $[0,\infty)$ such that
$a-1 \in C_c^\infty((0,R)).$ For $x \ne p,$ we let
$r := |x-p|, \ e := \frac{x-p}{|x-p|},$
and define
\begin{equation}
\label{eq:planarA}
A(x) := a(r)e\otimes e
  +\frac{1}{a(r)}\left(I-e\otimes e\right),
  \qquad A(p) := I.
\end{equation}
Then $A$ is smooth, symmetric, and uniformly elliptic,
with $\det A = 1.$ If
\begin{equation}
\label{eq:plancrit}
\int_0^R \frac{r}{a(r)} \, dr > \frac{3R^2}{2},
\end{equation}
then $L := \div(A\nabla\cdot)$ has weak mean value
balls which are not strong, as in
Proposition\refthm{radcrit}\!\!.

In particular, the condition holds for every such
profile satisfying
$$\frac{1}{4} \le a(r) \le 1,
  \qquad a(r) = \frac{1}{4}
  \quad \text{for } \frac{R}{8} \le r \le \frac{7R}{8}.$$
For these examples all eigenvalues of $A$ lie in
$[1/4,4].$
\end{corollary}

\begin{pf}
Define
$$T(r) := R\exp\left(-\int_r^R
                  \frac{dt}{t\,a(t)}\right)
  \quad (0 < r \le R),
  \qquad T(0) := 0,$$
and let $T(r) := r$ for $r \ge R.$ Then
$$\frac{T'(r)}{T(r)} = \frac{1}{r\,a(r)}.$$
The function $T$ is strictly increasing, equals a positive
constant times $r$ near zero, and equals $r$ near and beyond
$R.$ Therefore, with $s := T^{-1},$ the map
$$\Psi(p+te) := p+s(t)e,
  \qquad \Psi(p) := p,$$
is an orientation-preserving smooth radial diffeomorphism
equal to the identity outside a compact subset of $B_R(p).$
Under the conductivity transformation, the radial and
tangential eigenvalues are
$$\frac{s'(t)t}{s(t)} = a(s(t))
  \quad \text{and} \quad
  \frac{s(t)}{t s'(t)} = \frac{1}{a(s(t))},$$
respectively, so it is precisely $A.$
Proposition\refthm{radcrit}with $n = 2$ applies.
Finally, for the stated profiles,
$$\begin{aligned}
  \int_0^R \frac{r}{a(r)} \, dr
    &\ge \frac{R^2}{2}
       +3\int_{R/8}^{7R/8}r \, dr\\
    &= \frac{13R^2}{8} > \frac{3R^2}{2}
\end{aligned}$$
which verifies Equation\refeqn{plancrit}\!\!.
The eigenvalue bounds follow from Equation\refeqn{planarA}\!\!.
\end{pf}

\begin{remark}[Families of examples]
\label{radfamily}
The planar condition is strict and depends continuously
on $a$ in the uniform norm when the profiles are bounded
away from zero. It therefore persists under sufficiently
small uniform perturbations within the class of positive
smooth profiles with $a-1 \in C_c^\infty((0,R)).$
In particular, the sufficient condition holds on a nonempty
relatively open subset of this infinite-dimensional profile
class. This gives a precise sense in which the examples
are not isolated within the radial family; it does not
assert openness or genericity among all coefficient matrices.
More generally, Lemma~\refthm{radexist}shows directly that the
criterion in Proposition~\refthm{radcrit}has examples in every
dimension $n \ge 2.$
\end{remark}

\section{Enclosure and Holes in Canonical Mean Value Sets}
\label{holes}

\subsection{The enclosure principle}

In \cite{Ar}, the first author constructed a nonconvex mean value
set.  Here we prove an enclosure principle for weak mean value
sets which are not strong. For such a set $D,$ let $E$ be the
canonical set with the same center and volume. By
Corollary~\refthm{roughenc}\!\!, $\partial D \subset \overline{E}$
and $D$ contains a bounded component of the open exterior
$\R^n \setminus \overline{E}.$ 
With $\partial D$ of class $C^{1,1}$ and $A$ smooth near
$\partial D,$ Theorem~\refthm{ESMVWH}rules out contact between
$\partial D$ and $\partial E,$ so $\partial D \subset E$ and
$D$ contains a bounded component of the closed complement $E^c$
of positive measure.
Applied to the construction in the previous section,
this gives canonical sets with holes of positive
measure in every dimension $n \ge 2.$

Related hole formation occurs for weighted planar
Hele--Shaw flows. Ross and Witt Nystr\"{o}m
\cite[Proposition~1.4]{RWN2015} construct a smooth
positive area form on $\mathbb{P}^1$ whose associated
Hele--Shaw domains are multiply connected throughout
an interval of times.
Theorem~\refthm{ESMVWH}concerns Lebesgue averages for
divergence-form operators and derives the hole from
the failure of the strong mean value property.

The valid maximum-principle comparison from
\cite{ABthepain}, supplemented by a boundary-point
argument, yields the enclosure of the weak set's
boundary by the canonical set.
The ball example exhibits a configuration that the
flawed uniqueness proof incorrectly excluded.
Shapiro explicitly identified the failure of the
relevant boundary-point inference when the comparison
set has a hole; see \cite[pp.~33--34]{Sh}.
Equality of volumes and the measure-zero free boundary
then force a bounded complementary component of
positive measure.
For the Laplacian, an analogous comparison of the
potentials of harmonic and subharmonic balls appears
in \cite[Proposition 3.2]{ShSj}.

The boundary-point argument in Step 2 is related to
the Hopf-lemma comparisons at regular free boundary
points in \cite[proof of Lemma 4.1]{BBL} and
\cite[Section 1]{ArmB}.
Here an auxiliary elliptic problem on an annulus
inside an interior tangent ball supplies a linear
lower bound for the canonical potential near a
hypothetical common boundary point.
This contradicts the quadratic upper bound from
\cite[Theorem 3.2]{BH1}, without requiring regularity
of the canonical free boundary.

We now return from the specially constructed operators
arising from conductivity
transforms of Section\refthm{sharpness}to the general setting
of Section\refthm{AND}\!\!, with $L := \div( A \nabla \cdot )$
on $\R^n,$ where $A$ is bounded, measurable, symmetric, and
uniformly elliptic.

For the proof of Theorem~\refthm{ESMVWH}\!\!\!, write
$\mathcal{V} := |D|$ and $E := D_{\mathcal{V}^{1/n}}(x_0).$
Notice that no topological assumption on $D$ beyond the
connectedness forced by the weak mean value property is needed;
in particular, $D$ need not be simply connected and
$\R^n \setminus D$ need not be connected.

\begin{pf}
\textit{Step 1: Comparison of the potentials.}

Choose a bounded smooth ambient domain $\Omega$
compactly containing $D$ and $E.$
Let
$$\Upsilon_{D,E}
  := W_D-W_{E}$$
be the canonical comparison from
Lemma~\refthm{cancomp}\!\!.
That lemma gives
$$\Upsilon_{D,E} \le 0
  \qquad \text{in } \Omega.$$
Moreover,
$$L\Upsilon_{D,E}
  = \frac{1}{\mathcal{V}}
      \left(1-\bigrchi_{E}\right) \ge 0
  \qquad \text{in } D.$$
If $|D \triangle E| = 0,$ then Corollary\refthm{maxpotential}gives
$W_D = W_E \ge 0$ away from $x_0.$ Since $D$ is weak,
Theorem\refthm{openreps}would then imply that $D$ is strong.
Thus $|D \triangle E| > 0,$ and the strict comparison in
Corollary\refthm{maxpotential}gives
$$\Upsilon_{D,E} < 0 \qquad \text{throughout } D \cup E.$$

Figure~\ref{fig:hole} illustrates the enclosure geometry
and the volume balance that forces a hole.

\begin{figure}[t]
\centering
\includegraphics[width=0.75\textwidth]{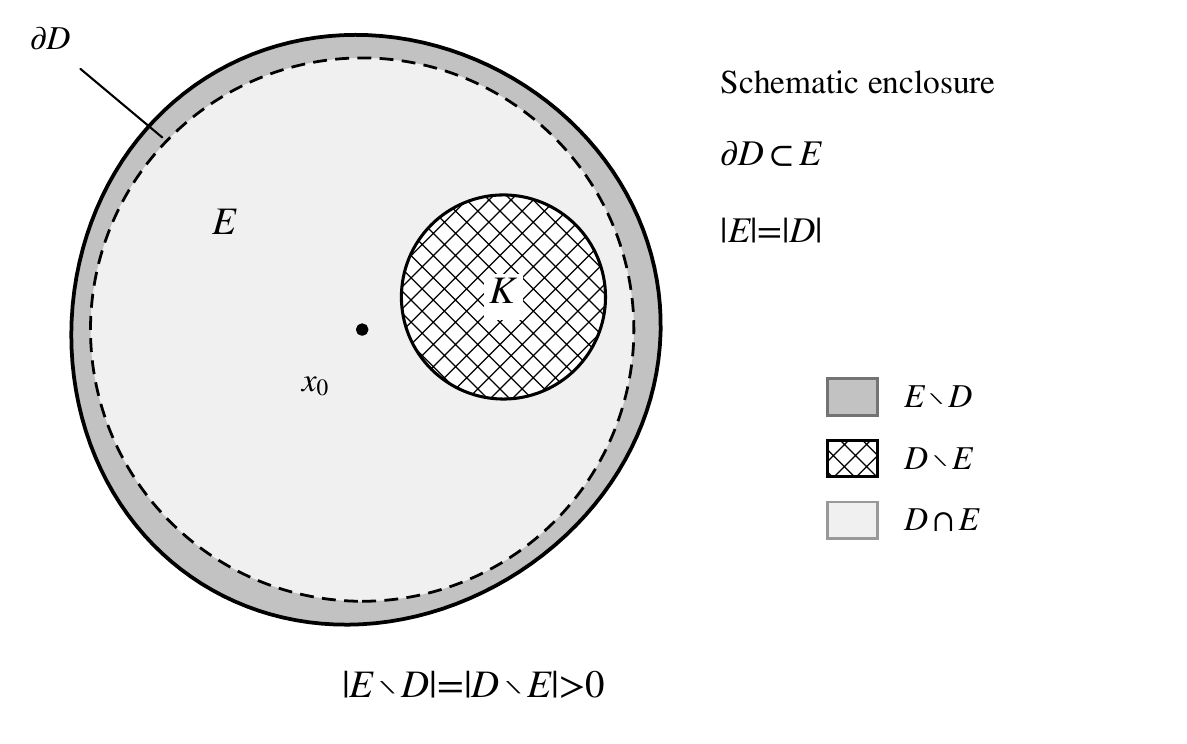}
\caption{Schematic enclosure geometry for the ball example.
The dashed curve is $\partial D,$ which lies entirely
inside $E.$
The region $E \setminus D$ has positive measure;
equality of the volumes gives
$|D \setminus E| = |E \setminus D| > 0.$
The component $K \subset D$ represents a bounded
complementary component of positive measure.
The drawing illustrates the set relations, not the
computed boundary of a canonical set or its number of holes.}
\label{fig:hole}
\end{figure}

\medskip
\noindent
\textit{Step 2: The canonical set contains $\partial D$.}

Choose an open neighborhood $N$ of $\partial D$ such that
$$\overline{N} \subset \Omega \setminus \{x_0\}$$
and $A$ is smooth on $N.$
Since $x_0 \notin N,$ Equation~\eqref{eq:OPHF} gives
$$LW_D
  = \frac{1}{\mathcal{V}}\bigrchi_D
  \qquad \text{in } N.$$
Interior elliptic regularity gives
$$W_D \in C^{1,\alpha}_{\mathrm{loc}}(N)
  \qquad \text{for every } 0 < \alpha < 1;$$
see \cite[Corollary 8.36]{GT}.
Since $W_D = 0$ outside $D,$ this regularity, together with the
$C^{1,1}$ boundary of $D,$ implies
$$W_D(z) = 0,
  \qquad \nabla W_D(z) = 0
  \qquad \text{for every } z \in \partial D.$$

Fix $z \in \partial D.$
Choose an interior tangent ball $B_R(q) \subset D$
with
$$z \in \partial B_R(q),
  \qquad \overline{B_R(q)} \subset N,$$
and set
$$\nu := \frac{z-q}{R},
  \qquad
  \mathcal{A}
  := B_R(q) \setminus \overline{B_{R/2}(q)}.$$
Let $H$ solve the smooth Dirichlet problem
$$
\begin{aligned}
\div(A\nabla H) &= 0
  &&\text{in } \mathcal{A},\\
H &= -1
  &&\text{on } \partial B_{R/2}(q),\\
H &= 0
  &&\text{on } \partial B_R(q).
\end{aligned}
$$
Then $H < 0$ in $\mathcal{A}.$
The boundary-point lemma applied to $H$ gives a
constant $c_0 > 0$ such that
$$H(z-t\nu) \le -c_0t$$
for all sufficiently small $t > 0;$
see \cite[Lemma 3.4]{GT}.

Set
$$\varepsilon
  := -\max_{\partial B_{R/2}(q)}
       \Upsilon_{D,E} > 0.$$
On the inner sphere,
$\Upsilon_{D,E} \le -\varepsilon,$
and on the outer sphere, $\Upsilon_{D,E} \le 0.$
Since $L\Upsilon_{D,E} \ge 0$ in $\mathcal{A},$
the weak comparison principle gives
$$\Upsilon_{D,E} \le \varepsilon H
  \qquad \text{in } \mathcal{A}.$$
Hence, with $c := \varepsilon c_0 > 0,$
$$\Upsilon_{D,E}(z-t\nu) \le -ct$$
for all sufficiently small $t > 0.$

Suppose now that $W_{E}(z) = 0.$
The vanishing of $W_D$ and its gradient at $z$ gives
$$W_D(z-t\nu) = o(t).$$
It follows that
$$
\begin{aligned}
W_{E}(z-t\nu)
 &= W_D(z-t\nu)-\Upsilon_{D,E}(z-t\nu)\\
  &\ge ct+o(t).
\end{aligned}
$$
Thus $W_{E}(z-t\nu) > 0$ for all sufficiently
small $t > 0,$ and consequently $z \in \partial E.$
Because $z \ne x_0,$ the quadratic upper bound from
\cite[Theorem 3.2]{BH1} applies and gives
$$W_{E}(z-t\nu) \le Ct^2$$
for sufficiently small $t > 0.$
This contradicts the preceding linear lower bound.
Therefore,
$$\partial D \subset E.$$

\medskip
\noindent
\textit{Step 3: A bounded complementary component
has positive measure.} \nopagebreak[4]

Since $|D| = |E|$ and $|\partial E| = 0,$ Step~1 gives
$$|D \setminus \overline{E}| = |D \setminus E|
  = \frac{1}{2}|D \triangle E| > 0.$$
The open set $D \setminus \overline{E}$ therefore contains a ball.
Let $K$ be the connected component of $E^c$
containing that ball.
The set $K$ meets $D$ and avoids $\partial D,$
because $\partial D \subset E.$
Connectedness therefore implies $K \subset D.$
Thus $K$ is bounded and has positive measure.

Finally, $E$ is bounded, so its complement contains
the exterior of a sufficiently large ball.
Since $n \ge 2,$ that exterior is connected and has
infinite measure, and therefore lies in an unbounded
component of $E^c$ of infinite measure.

For the operator in Theorem~\refthm{WDNIS}\!\!,
the ball $D = B_1(0)$ is weak but not strong, and
$A = I$ in a full neighborhood of $\partial D.$
Thus the hypotheses of the theorem are satisfied,
and the final assertion follows.
\end{pf}

\begin{remark}[One-sided coefficient regularity]
\label{onesidedregularity}
The smoothness assumption on $A$ near $\partial D$ in
Theorem~\refthm{ESMVWH}is made for convenience and is
stronger than necessary. More precisely, the same argument applies
if, for some $0 < \alpha < 1,$ there exist an open neighborhood
$N$ of $\partial D$ and a symmetric, uniformly elliptic matrix
field
$$\widetilde A \in C^{0,\alpha}(N)$$
such that
$$A = \widetilde A
  \qquad \text{almost everywhere in } N \cap D.$$
Indeed, since $W_D = 0$ outside $D,$ we have
$\nabla W_D = 0$ almost everywhere there, and therefore
$$\widetilde A\nabla W_D = A\nabla W_D
  \qquad \text{almost everywhere in } N.$$
After shrinking $N$ to avoid $x_0,$ it follows that
$$\div(\widetilde A\nabla W_D)
  = \frac{1}{\mathcal{V}}\bigrchi_D
  \qquad \text{in } N.$$
The required $C^{1,\alpha}$ regularity then follows from
\cite[Corollary 8.36]{GT}, and the boundary-point estimate used
in Step 2 remains valid for divergence-form equations with
\Holder continuous leading coefficients. In dimension two this
goes back to Finn and Gilbarg \cite{FG57}; for a
general-dimensional result under the weaker Dini-continuity
hypothesis, see Apushkinskaya and Nazarov \cite{AN19}.
Thus the actual coefficients outside $D$ may change abruptly
across $\partial D$ and need only satisfy the standing
coefficient assumptions.  In particular, for the ball in
Theorem~\refthm{WDNIS}\!\!\!, the conclusion remains valid after
any bounded measurable, symmetric, uniformly elliptic
modification of $A$ outside $D.$
\end{remark}

\begin{remark}[Localization of the interior holes]
\label{holelocation}
Under the hypotheses of Theorem~\refthm{ESMVWH}\!\!,
$D$ is weak but not strong, so Theorem~\refthm{WisSiffPisP}and
Corollary~\ref{maxpotential} give
$$D \setminus E
  \subset \{x \in D : W_D(x) < 0\}.$$

For the unit-ball counterexample obtained in the proof of
Theorem~\refthm{WDNIS}\!\!, let $D_* := B_\varrho(x_0)$ be the
pre-rescaling ball supplied by Proposition~\refthm{radcrit}\!\!,
and let $T(z) := x_0 + \varrho z$ be the affine rescaling used in
the proof of Theorem~\refthm{WDNIS}\!\!. With the notation from
the proof of Proposition~\refthm{radcrit}\!\!, set
$$\widehat p := T^{-1}(p),
  \qquad \widehat R := \frac{R}{\varrho},$$
and
$$\widehat P(z) := \frac{1}{\varrho^2}P_{D_*}(T(z)),
  \qquad
  \widehat\phi(z) := \frac{1}{\varrho^2}\phi(T(z)),
  \qquad
  \widehat Q := \widehat P + \widehat\phi.$$
After the renaming in the proof of Theorem~\refthm{WDNIS}\!\!,
the rescaling gives
$$W_D
  = \frac{\widehat Q \circ \Psi^{-1}}{\omega_n}.$$
Moreover, $\widehat P > 0$ in $B_1(0) \setminus \{0\},$
$\widehat\phi = 0$ outside
$B_{\widehat R}(\widehat p),$ and $\Psi$ preserves
$B_{\widehat R}(\widehat p).$ Hence
$$D \setminus E
  \subset B_{\widehat R}(\widehat p).$$
Thus every connected component of $E^c$ contained in $D$ lies
in the rescaled perturbation region.
\end{remark}

\subsection{Application to anisotropic inverse conductivity}

The preceding construction also gives a consequence for
the anisotropic inverse conductivity problem. For a smooth
conductivity matrix $\sigma$ on a bounded smooth domain
$\Omega,$ write $\Lambda_\sigma^\Omega$ for the
Dirichlet-to-Neumann map
$$\Lambda_\sigma^\Omega f
 := \left. \nu \cdot \sigma\nabla u_f\right|_{\partial\Omega}, \ \ \
  \div(\sigma\nabla u_f) = 0 \ \text{in }\Omega, \ \ \
  u_f = f\ \text{on }\partial \Omega,$$
where $\nu$ is the outward unit normal.

\begin{corollary}[Boundary measurements and mean value set topology]
\label{conductivity}
For every $n \ge 2,$ there exist a bounded smooth domain
$\Omega \subset \R^n$ containing $0$ and a smooth, symmetric,
uniformly elliptic matrix field $A$ on $\R^n$ such that
$A = I$ near $0$ and outside a compact subset of $\Omega,$ and
$$\Lambda_A^\Omega = \Lambda_I^\Omega.$$
Fix the center $0$ and the volume $\omega_n := |B_1(0)|.$
The complement of the canonical mean value set for
$\div(A\nabla\cdot)$ has a bounded component of positive measure.
The canonical mean value set for $\Delta$ is $B_1(0),$
whose complement is connected.
Both canonical sets are compactly contained in $\Omega.$
When $n = 2,$ one may additionally require $\det A = 1.$
All averages here are taken with respect to Lebesgue measure.
\end{corollary}

\begin{pf}
Take $A$ and $\Psi$ from the proof of
Theorem~\refthm{WDNIS}\!\!\!. Thus $A$ is the conductivity
pushforward of $I$ by $\Psi,$ and $\Psi$ is the
identity near $0$ and outside a compact subset of
$B_1(0).$
Choose a ball $\Omega$ compactly containing both
$B_1(0)$ and the canonical set
$E := D_{\omega_n^{1/n}}(0)$ for $\div(A\nabla\cdot).$
Since $\Psi$ is the identity near $\partial\Omega,$
invariance under changes of variables that fix the
boundary gives $\Lambda_A^\Omega = \Lambda_I^\Omega;$
see \cite[Section 2.2]{KSVW} and, in the plane,
\cite[Section 1]{ALP}.

Theorem~\refthm{ESMVWH}gives a bounded component of
$E^c$ of positive measure. For the Laplacian,
the canonical set with the same center and volume is
$B_1(0),$ whose complement is connected and unbounded.
The determinant assertion follows from
Theorem~\refthm{WDNIS}\!\!.
\end{pf}
The fixed choice of Lebesgue measure is essential:
under $x = \Psi(y),$ the measure $dx$ pulls back to
$\det D\Psi(y) \, dy.$ Thus the conductivity change of
variables preserves boundary measurements while
introducing a nonconstant density into the averaging
formula.

\bibliographystyle{plain}
\bibliography{BIB}

@article{Abat26,
  AUTHOR = {Abatangelo, Nicola},
  TITLE = {Rigidity of balls in the solid mean value property
           for polyharmonic functions},
  JOURNAL = {Bull. Lond. Math. Soc.},
  VOLUME = {58},
  YEAR = {2026},
  NUMBER = {6},
  PAGES = {e70416},
  DOI = {10.1112/blms.70416},
}

@article{AbbBon13,
  author  = {Abbondanza, Beatrice and Bonfiglioli, Andrea},
  title   = {The {Dirichlet} problem and the inverse mean-value
             theorem for a class of divergence form operators},
  journal = {J. Lond. Math. Soc. (2)},
  volume  = {87},
  number  = {2},
  pages   = {321--346},
  year    = {2013},
  doi     = {10.1112/jlms/jds050}
}

@book{AMC,
  author    = {Ambrosio, Luigi and Carlotto, Alessandro
               and Massaccesi, Annalisa},
  title     = {Lectures on Elliptic Partial Differential Equations},
  series    = {Lecture Notes. Scuola Normale Superiore di Pisa
               (New Series)},
  volume    = {18},
  publisher = {Edizioni della Normale},
  address   = {Pisa},
  year      = {2018},
  doi       = {10.1007/978-88-7642-651-3}
}

@article{AHT,
  author = {Ameur, Yacin and Helmer, Martin and Tellander, Felix},
  title = {On the Uniqueness Problem for Quadrature Domains},
  journal = {Computational Methods and Function Theory},
  volume = {21},
  pages = {473--504},
  year = {2021},
  doi = {10.1007/s40315-021-00373-w}
}

@article{AN19,
  AUTHOR = {Apushkinskaya, Darya E. and Nazarov, Alexander I.},
  TITLE = {On the boundary point principle for divergence-type equations},
  JOURNAL = {Atti Accad. Naz. Lincei Rend. Lincei Mat. Appl.},
  VOLUME = {30},
  YEAR = {2019},
  NUMBER = {4},
  PAGES = {677--699},
  DOI = {10.4171/RLM/867},
}

@article{ArmB,
    AUTHOR = {Armstrong, Niles and Blank, Ivan},
     TITLE = {Nondegenerate motion of singular points in obstacle problems
              with varying data},
   JOURNAL = {J. Differential Equations},
  FJOURNAL = {Journal of Differential Equations},
    VOLUME = {267},
      YEAR = {2019},
    NUMBER = {9},
     PAGES = {5370--5387},
      ISSN = {0022-0396},
   MRCLASS = {Prelim},
  MRNUMBER = {3991562},
       DOI = {10.1016/j.jde.2019.05.032},
       URL = {https://doi-org.er.lib.k-state.edu/10.1016/j.jde.2019.05.032},
}

@article{ABthepain,
    AUTHOR = {Armstrong, Niles and Blank, Ivan},
     TITLE = {A uniqueness theorem for mean value sets for elliptic divergence form operators},
   JOURNAL = {arXiv:1907.12523},
      YEAR = {2019},
}

@unpublished{ABpartII,
  author = {Niles Armstrong and Ivan Blank},
  title = {Formation, Evolution, and Extinction of Holes in
           Canonical Mean Value Sets},
  note = {Companion manuscript},
  year = {2026}
}

@article {Ar,
    AUTHOR = {Armstrong, Niles},
     TITLE = {Properties of mean value sets: angle conditions, blowup
              solutions, and nonconvexity},
   JOURNAL = {Potential Anal.},
  FJOURNAL = {Potential Analysis. An International Journal Devoted to the
              Interactions between Potential Theory, Probability Theory,
              Geometry and Functional Analysis},
    VOLUME = {52},
      YEAR = {2020},
    NUMBER = {3},
     PAGES = {527--544},
      ISSN = {0926-2601,1572-929X},
   MRCLASS = {35J15 (35B05 35B44 35J08 35R05 35R35)},
  MRNUMBER = {4067302},
MRREVIEWER = {Roberto\ Ognibene},
       DOI = {10.1007/s11118-018-9741-3},
       URL = {https://doi-org.er.lib.k-state.edu/10.1007/s11118-018-9741-3},
}

@article {AryB,
    AUTHOR = {Aryal, Ashok and Blank, Ivan},
     TITLE = {Geometry of mean value sets for general divergence form
              uniformly elliptic operators},
   JOURNAL = {Potential Anal.},
  FJOURNAL = {Potential Analysis. An International Journal Devoted to the
              Interactions between Potential Theory, Probability Theory,
              Geometry and Functional Analysis},
    VOLUME = {50},
      YEAR = {2019},
    NUMBER = {1},
     PAGES = {43--54},
      ISSN = {0926-2601},
   MRCLASS = {35J86 (35J08 35J15 35J20 35R35)},
  MRNUMBER = {3900845},
       DOI = {10.1007/s11118-017-9672-4},
       URL = {https://doi-org.er.lib.k-state.edu/10.1007/s11118-017-9672-4},
}

@article {ALP,
    AUTHOR = {Astala, Kari and P\"aiv\"arinta, Lassi and Lassas, Matti},
     TITLE = {Calder\'on's inverse problem for anisotropic conductivity in
              the plane},
   JOURNAL = {Comm. Partial Differential Equations},
  FJOURNAL = {Communications in Partial Differential Equations},
    VOLUME = {30},
      YEAR = {2005},
    NUMBER = {1-3},
     PAGES = {207--224},
      ISSN = {0360-5302,1532-4133},
   MRCLASS = {35R30 (35J25)},
  MRNUMBER = {2131051},
MRREVIEWER = {Yongzhi\ Xu},
       DOI = {10.1081/PDE-200044485},
       URL = {https://doi-org.er.lib.k-state.edu/10.1081/PDE-200044485},
}

@article {BBL,
AUTHOR = {Benson, Brian and Blank, Ivan and LeCrone, Jeremy},
     TITLE = {Mean value theorems for {R}iemannian manifolds via the obstacle problem},
   JOURNAL = {J. Geom. Anal.},
  FJOURNAL = {Journal of Geometric Analysis},
    VOLUME = {29},
      YEAR = {2019},
    NUMBER = {3},
     PAGES = {2752--2775},
      ISSN = {1050-6926},
   MRCLASS = {58J05 (35B05 35J15 35R35)},
  MRNUMBER = {3969443},
       DOI = {10.1007/s12220-018-0093-4},
       URL = {https://doi-org.er.lib.k-state.edu/10.1007/s12220-018-0093-4},
}

@article {BH1,
    AUTHOR = {Blank, Ivan and Hao, Zheng},
     TITLE = {The mean value theorem and basic properties of the obstacle
              problem for divergence form elliptic operators},
   JOURNAL = {Comm. Anal. Geom.},
  FJOURNAL = {Communications in Analysis and Geometry},
    VOLUME = {23},
      YEAR = {2015},
    NUMBER = {1},
     PAGES = {129--158},
      ISSN = {1019-8385},
   MRCLASS = {35J15 (35B05 35J08)},
  MRNUMBER = {3291366},
MRREVIEWER = {Juraj F{\"o}ldes},
       DOI = {10.4310/CAG.2015.v23.n1.a4},
       URL = {http://dx.doi.org/10.4310/CAG.2015.v23.n1.a4},
}

@article {BH2,
    AUTHOR = {Blank, Ivan and Hao, Zheng},
     TITLE = {Reifenberg flatness of free boundaries in obstacle problems
              with {VMO} ingredients},
   JOURNAL = {Calc. Var. Partial Differential Equations},
  FJOURNAL = {Calculus of Variations and Partial Differential Equations},
    VOLUME = {53},
      YEAR = {2015},
    NUMBER = {3-4},
     PAGES = {943--959},
      ISSN = {0944-2669},
   MRCLASS = {35J87 (35A01 35A02 35B65 35J15 35J20)},
  MRNUMBER = {3347485},
MRREVIEWER = {Rodica Luca},
       DOI = {10.1007/s00526-014-0772-3},
       URL = {https://doi-org.er.lib.k-state.edu/10.1007/s00526-014-0772-3},
}

@article {B,
    AUTHOR = {Blank, Ivan},
     TITLE = {Sharp results for the regularity and stability of the free
              boundary in the obstacle problem},
   JOURNAL = {Indiana Univ. Math. J.},
  FJOURNAL = {Indiana University Mathematics Journal},
    VOLUME = {50},
      YEAR = {2001},
    NUMBER = {3},
     PAGES = {1077--1112},
      ISSN = {0022-2518},
     CODEN = {IUMJAB},
   MRCLASS = {35R35 (49J40)},
  MRNUMBER = {1871348},
MRREVIEWER = {Martin Fuchs},
       DOI = {10.1512/iumj.2001.50.1906},
       URL = {http://dx.doi.org/10.1512/iumj.2001.50.1906},
}

@article{BRG25,
  AUTHOR = {Bou-Rabee, Ahmed and Gwynne, Ewain},
  TITLE = {Harmonic balls in {L}iouville quantum gravity},
  JOURNAL = {Proc. Lond. Math. Soc.},
  VOLUME = {130},
  YEAR = {2025},
  NUMBER = {1},
  PAGES = {e70018},
  DOI = {10.1112/plms.70018},
}

@article{BW91,
  author  = {Brummelhuis, Raymond and Wiegerinck, Jan},
  title   = {Representing measures for the disc algebra and for the ball algebra},
  journal = {Ann. Polon. Math.},
  volume  = {55},
  number  = {1},
  pages   = {19--35},
  year    = {1991},
  doi     = {10.4064/ap-55-1-19-35}
}

@book {C2,
    AUTHOR = {Caffarelli, Luis A.},
     TITLE = {The obstacle problem},
    SERIES = {Lezioni Fermiane. [Fermi Lectures]},
 PUBLISHER = {Accademia Nazionale dei Lincei, Rome; Scuola Normale
              Superiore, Pisa},
      YEAR = {1998},
     PAGES = {ii+54},
   MRCLASS = {49-02 (35J20 35J85 35R35 49J40 49J45 49K20)},
  MRNUMBER = {2011808},
MRREVIEWER = {Martin Fuchs},
}

@article {CR,
    AUTHOR = {Caffarelli, Luis A. and Roquejoffre, Jean-Michel},
     TITLE = {Uniform {H}\"{o}lder estimates in a class of elliptic systems and
              applications to singular limits in models for diffusion
              flames},
   JOURNAL = {Arch. Ration. Mech. Anal.},
  FJOURNAL = {Archive for Rational Mechanics and Analysis},
    VOLUME = {183},
      YEAR = {2007},
    NUMBER = {3},
     PAGES = {457--487},
      ISSN = {0003-9527},
   MRCLASS = {35J55 (35B65 80A25)},
  MRNUMBER = {2278412},
MRREVIEWER = {Helmut Abels},
       DOI = {10.1007/s00205-006-0013-9},
       URL = {https://doi-org.er.lib.k-state.edu/10.1007/s00205-006-0013-9},
}

@article{CL25,
  AUTHOR = {Cupini, Giovanni and Lanconelli, Ermanno},
  TITLE = {On the harmonic characterization of the spheres:
           a sharp stability inequality},
  JOURNAL = {Potential Anal.},
  VOLUME = {63},
  YEAR = {2025},
  PAGES = {1181--1207},
  DOI = {10.1007/s11118-025-10200-9},
}

@article{CFLZ,
  author  = {Cupini, Giovanni and Lanconelli, Ermanno and
             Fusco, Nicola and Zhong, Xiao},
  title   = {A sharp stability result for the {Gauss}
             mean value formula},
  journal = {J. Anal. Math.},
  volume  = {142},
  number  = {2},
  year    = {2020},
  pages   = {587--603},
  doi     = {10.1007/s11854-020-0144-1},
  url     = {https://doi.org/10.1007/s11854-020-0144-1}
}

@article{Dal08,
  author  = {Dalmasso, Robert},
  title   = {On the mean-value property of superharmonic functions},
  journal = {Ann. Polon. Math.},
  volume  = {93},
  number  = {1},
  pages   = {85--90},
  year    = {2008},
  doi     = {10.4064/ap93-1-7}
}

@article {EpSchiff,
    AUTHOR = {Epstein, Bernard and Schiffer, M. M.},
     TITLE = {On the mean-value property of harmonic functions},
   JOURNAL = {J. Analyse Math.},
  FJOURNAL = {Journal d'Analyse Math\'{e}matique},
    VOLUME = {14},
      YEAR = {1965},
     PAGES = {109--111},
      ISSN = {0021-7670},
   MRCLASS = {31.11},
  MRNUMBER = {0177124},
MRREVIEWER = {John W. Green},
       DOI = {10.1007/BF02806381},
       URL = {https://doi-org.er.lib.k-state.edu/10.1007/BF02806381},
}

@article {Ep,
    AUTHOR = {Epstein, Bernard},
     TITLE = {On the mean-value property of harmonic functions},
   JOURNAL = {Proc. Amer. Math. Soc.},
  FJOURNAL = {Proceedings of the American Mathematical Society},
    VOLUME = {13},
      YEAR = {1962},
     PAGES = {830},
      ISSN = {0002-9939},
   MRCLASS = {31.10},
  MRNUMBER = {140700},
MRREVIEWER = {John W. Green},
       DOI = {10.2307/2034188},
       URL = {https://doi-org.er.lib.k-state.edu/10.2307/2034188},
}

@article{FG57,
  AUTHOR = {Finn, R. and Gilbarg, D.},
  TITLE = {Asymptotic behavior and uniqueness of plane subsonic flows},
  JOURNAL = {Comm. Pure Appl. Math.},
  VOLUME = {10},
  YEAR = {1957},
  PAGES = {23--63},
  DOI = {10.1002/cpa.3160100102},
}

@article {GardSj,
    AUTHOR = {Gardiner, Stephen J. and Sj\"odin, Tomas},
     TITLE = {Partial balayage for the {H}elmholtz equation},
   JOURNAL = {Potential Anal.},
  FJOURNAL = {Potential Analysis. An International Journal Devoted to the
              Interactions between Potential Theory, Probability Theory,
              Geometry and Functional Analysis},
    VOLUME = {63},
      YEAR = {2025},
    NUMBER = {4},
     PAGES = {1671--1697},
      ISSN = {0926-2601,1572-929X},
   MRCLASS = {35J05 (31B15)},
  MRNUMBER = {4990500},
       DOI = {10.1007/s11118-025-10217-0},
       URL = {https://doi-org.er.lib.k-state.edu/10.1007/s11118-025-10217-0},
}

@article {GO,
    AUTHOR = {Goldstein, Myron and Ow, Wellington H.},
     TITLE = {On the mean-value property of harmonic functions},
   JOURNAL = {Proc. Amer. Math. Soc.},
  FJOURNAL = {Proceedings of the American Mathematical Society},
    VOLUME = {29},
      YEAR = {1971},
     PAGES = {341--344},
      ISSN = {0002-9939},
   MRCLASS = {31.10},
  MRNUMBER = {279320},
MRREVIEWER = {M. M. Schiffer},
       DOI = {10.2307/2038138},
       URL = {https://doi-org.er.lib.k-state.edu/10.2307/2038138},
}

@book {GT,
    AUTHOR = {Gilbarg, David and Trudinger, Neil S.},
     TITLE = {Elliptic partial differential equations of second order},
    SERIES = {Classics in Mathematics},
      NOTE = {Reprint of the 1998 edition},
 PUBLISHER = {Springer-Verlag, Berlin},
      YEAR = {2001},
     PAGES = {xiv+517},
      ISBN = {3-540-41160-7},
   MRCLASS = {35-02 (35Jxx)},
  MRNUMBER = {1814364},
}

@article{Gus88,
  author  = {Gustafsson, Bj{\"o}rn},
  title   = {Singular and special points on quadrature domains
             from an algebraic geometric point of view},
  journal = {J. Anal. Math.},
  volume  = {51},
  number  = {1},
  pages   = {91--117},
  year    = {1988},
  doi     = {10.1007/BF02791120}
}

@article {GRooooos,
    AUTHOR = {Gustafsson, Bj\"orn and Roos, Joakim},
     TITLE = {Partial balayage on {R}iemannian manifolds},
   JOURNAL = {J. Math. Pures Appl. (9)},
  FJOURNAL = {Journal de Math\'ematiques Pures et Appliqu\'ees. Neuvi\`eme
              S\'erie},
    VOLUME = {118},
      YEAR = {2018},
     PAGES = {82--127},
      ISSN = {0021-7824,1776-3371},
   MRCLASS = {31C12 (35R01 35R35 58A14)},
  MRNUMBER = {3852470},
MRREVIEWER = {Mohamed\ Selmi},
       DOI = {10.1016/j.matpur.2017.07.013},
       URL = {https://doi-org.er.lib.k-state.edu/10.1016/j.matpur.2017.07.013},
}

@incollection {GusShap,
    AUTHOR = {Gustafsson, Bj\"orn and Shapiro, Harold S.},
     TITLE = {What is a quadrature domain?},
 BOOKTITLE = {Quadrature domains and their applications},
    SERIES = {Oper. Theory Adv. Appl.},
    VOLUME = {156},
     PAGES = {1--25},
 PUBLISHER = {Birkh\"auser, Basel},
      YEAR = {2005},
      ISBN = {3-7643-7145-5},
   MRCLASS = {30-02 (31-02 35A21 35R35 47B20 76M40)},
  MRNUMBER = {2129734},
MRREVIEWER = {Darren\ G.\ Crowdy},
       DOI = {10.1007/3-7643-7316-4\_1},
       URL = {https://doi-org.er.lib.k-state.edu/10.1007/3-7643-7316-4_1},
}

@book {HL,
    AUTHOR = {Han, Qing and Lin, Fanghua},
     TITLE = {Elliptic partial differential equations},
    SERIES = {Courant Lecture Notes in Mathematics},
    VOLUME = {1},
   EDITION = {Second},
 PUBLISHER = {Courant Institute of Mathematical Sciences, New York; American
              Mathematical Society, Providence, RI},
      YEAR = {2011},
     PAGES = {x+147},
      ISBN = {978-0-8218-5313-9},
   MRCLASS = {35Jxx (35-01 35B50)},
  MRNUMBER = {2777537},
}

@article{HJS02,
  author  = {Hedenmalm, H{\aa}kan and Jakobsson, Stefan and Shimorin, Sergei},
  title   = {A biharmonic maximum principle for hyperbolic surfaces},
  journal = {J. Reine Angew. Math.},
  volume  = {550},
  pages   = {25--75},
  year    = {2002},
  doi     = {10.1515/crll.2002.074}
}

@article{HS02,
  author  = {Hedenmalm, H{\aa}kan and Shimorin, Sergei},
  title   = {{Hele--Shaw} flow on hyperbolic surfaces},
  journal = {J. Math. Pures Appl. (9)},
  volume  = {81},
  number  = {3},
  pages   = {187--222},
  year    = {2002},
  doi     = {10.1016/S0021-7824(01)01222-3}
}

@book {KS,
    AUTHOR = {Kinderlehrer, David and Stampacchia, Guido},
     TITLE = {An introduction to variational inequalities and their
              applications},
    SERIES = {Pure and Applied Mathematics},
    VOLUME = {88},
 PUBLISHER = {Academic Press, Inc. [Harcourt Brace Jovanovich, Publishers],
              New York-London},
      YEAR = {1980},
     PAGES = {xiv+313},
      ISBN = {0-12-407350-6},
   MRCLASS = {49A29 (35R35 47H15)},
  MRNUMBER = {567696},
MRREVIEWER = {Erich Miersemann},
}

@article{KSVW,
  AUTHOR = {Kohn, R. V. and Shen, H. and
            Vogelius, M. S. and Weinstein, M. I.},
  TITLE = {Cloaking via change of variables in electric
           impedance tomography},
  JOURNAL = {Inverse Problems},
  VOLUME = {24},
  YEAR = {2008},
  NUMBER = {1},
  PAGES = {015016, 21},
  DOI = {10.1088/0266-5611/24/1/015016}
}

@article{Koosis83,
  author  = {Koosis, Paul},
  title   = {La plus petite majorante surharmonique et son rapport
             avec l'existence des fonctions enti{\`e}res de type
             exponentiel jouant le r{\^o}le de multiplicateurs},
  journal = {Ann. Inst. Fourier (Grenoble)},
  volume  = {33},
  number  = {1},
  pages   = {67--107},
  year    = {1983},
  doi     = {10.5802/aif.905}
}

@article{KSS26,
  AUTHOR = {Kow, Pu-Zhao and Shahgholian, Henrik and Sj{\"o}din, Tomas},
  TITLE = {Multiphase quadrature domains (existence and uniqueness)},
  JOURNAL = {Anal. Math. Phys.},
  VOLUME = {16},
  YEAR = {2026},
  NOTE = {Article 88},
  DOI = {10.1007/s13324-026-01232-4},
}

@article{KLSS24,
  AUTHOR = {Kow, Pu-Zhao and Larson, Simon and Salo, Mikko
            and Shahgholian, Henrik},
  TITLE = {Quadrature domains for the {H}elmholtz equation with
           applications to non-scattering phenomena},
  JOURNAL = {Potential Anal.},
  VOLUME = {60},
  YEAR = {2024},
  NUMBER = {1},
  PAGES = {387--424},
  DOI = {10.1007/s11118-022-10054-5},
}

@article {Ku,
    AUTHOR = {Kuran, \"U.},
     TITLE = {On the mean-value property of harmonic functions},
   JOURNAL = {Bull. London Math. Soc.},
  FJOURNAL = {The Bulletin of the London Mathematical Society},
    VOLUME = {4},
      YEAR = {1972},
     PAGES = {311--312},
      ISSN = {0024-6093},
   MRCLASS = {31B20},
  MRNUMBER = {0320348},
MRREVIEWER = {W. H. Ow},
       DOI = {10.1112/blms/4.3.311},
       URL = {http://dx.doi.org/10.1112/blms/4.3.311},
}

@article {LSW,
    AUTHOR = {Littman, W. and Stampacchia, G. and Weinberger, H. F.},
     TITLE = {Regular points for elliptic equations with discontinuous
              coefficients},
   JOURNAL = {Ann. Scuola Norm. Sup. Pisa (3)},
    VOLUME = {17},
      YEAR = {1963},
     PAGES = {43--77},
   MRCLASS = {35.42},
  MRNUMBER = {0161019},
MRREVIEWER = {C. B. Morrey, Jr.},
}

@article{Loutsenko06,
  author  = {Loutsenko, Igor},
  title   = {The variable coefficient {Hele-Shaw} problem,
             integrability and quadrature identities},
  journal = {Comm. Math. Phys.},
  volume  = {268},
  number  = {2},
  pages   = {465--479},
  year    = {2006},
  doi     = {10.1007/s00220-006-0099-9}
}

@article{MY26,
  AUTHOR = {Morav{\'i}k, Mat\v{e}j and Youssfi, El Hassan},
  TITLE = {A converse mean-value property},
  JOURNAL = {arXiv:2606.22199},
  YEAR = {2026},
}

@article {N,
    AUTHOR = {Ni, Lei},
     TITLE = {Mean value theorems on manifolds},
   JOURNAL = {Asian J. Math.},
  FJOURNAL = {Asian Journal of Mathematics},
    VOLUME = {11},
      YEAR = {2007},
    NUMBER = {2},
     PAGES = {277--304},
      ISSN = {1093-6106},
   MRCLASS = {58J35 (35J65 53C44)},
  MRNUMBER = {2328895},
MRREVIEWER = {Xi Ping Zhu},
       DOI = {10.4310/AJM.2007.v11.n2.a6},
       URL = {http://dx.doi.org.er.lib.k-state.edu/10.4310/AJM.2007.v11.n2.a6},
}

@article{PP24,
  AUTHOR = {Pallara, Diego and Polidoro, Sergio},
  TITLE = {Mean value formulas for classical solutions to some
           degenerate elliptic equations in {C}arnot groups},
  JOURNAL = {Discrete Contin. Dyn. Syst. Ser. S},
  VOLUME = {17},
  YEAR = {2024},
  NUMBER = {5-6},
  PAGES = {1733--1749},
  DOI = {10.3934/dcdss.2022144},
}

@article{RWN2015,
  author  = {Ross, Julius and Witt Nystr{\"o}m, David},
  title   = {Harmonic discs of solutions to the complex
             homogeneous {Monge--Amp\`ere} equation},
  journal = {Publ. Math. Inst. Hautes \'{E}tudes Sci.},
  volume  = {122},
  pages   = {315--335},
  year    = {2015},
  doi     = {10.1007/s10240-015-0074-0}
}

@book {S,
    AUTHOR = {Sakai, Makoto},
     TITLE = {Quadrature domains},
    SERIES = {Lecture Notes in Mathematics},
    VOLUME = {934},
 PUBLISHER = {Springer-Verlag, Berlin-New York},
      YEAR = {1982},
     PAGES = {i+133},
      ISBN = {3-540-11562-5},
   MRCLASS = {41A55 (31A05)},
  MRNUMBER = {663007},
MRREVIEWER = {M. Goldstein},
}

@article {Shah,
    AUTHOR = {Shahgholian, Henrik},
     TITLE = {On quadrature domains and the {S}chwarz potential},
   JOURNAL = {J. Math. Anal. Appl.},
  FJOURNAL = {Journal of Mathematical Analysis and Applications},
    VOLUME = {171},
      YEAR = {1992},
    NUMBER = {1},
     PAGES = {61--78},
      ISSN = {0022-247X},
   MRCLASS = {31B05 (30E99 31B20 46F10)},
  MRNUMBER = {1192493},
       DOI = {10.1016/0022-247X(92)90376-O},
       URL = {https://doi-org.er.lib.k-state.edu/10.1016/0022-247X(92)90376-O},
}

@article {ShSj,
    AUTHOR = {Shahgholian, Henrik and Sj\"odin, Tomas},
     TITLE = {Harmonic balls and the two-phase {S}chwarz function},
   JOURNAL = {Complex Var. Elliptic Equ.},
  FJOURNAL = {Complex Variables and Elliptic Equations. An International
              Journal},
    VOLUME = {58},
      YEAR = {2013},
    NUMBER = {6},
     PAGES = {837--852},
      ISSN = {1747-6933,1747-6941},
   MRCLASS = {31B05 (30E20 30E25 31B20 35R35)},
  MRNUMBER = {3170664},
MRREVIEWER = {Stephen\ J.\ Gardiner},
       DOI = {10.1080/17476933.2011.622046},
       URL = {https://doi-org.er.lib.k-state.edu/10.1080/17476933.2011.622046},
}

@book {Sh,
    AUTHOR = {Shapiro, Harold S.},
     TITLE = {The {S}chwarz function and its generalization to higher
              dimensions},
    SERIES = {University of Arkansas Lecture Notes in the Mathematical
              Sciences},
    VOLUME = {9},
      NOTE = {A Wiley-Interscience Publication},
 PUBLISHER = {John Wiley \& Sons, Inc., New York},
      YEAR = {1992},
     PAGES = {xiv+108},
      ISBN = {0-471-57127-X},
   MRCLASS = {30E99 (30C40 30E25 46E20)},
  MRNUMBER = {1160990},
MRREVIEWER = {Steven George Krantz},
}

@article {W,
    AUTHOR = {Weiss, Georg S.},
     TITLE = {A homogeneity improvement approach to the obstacle problem},
   JOURNAL = {Invent. Math.},
  FJOURNAL = {Inventiones Mathematicae},
    VOLUME = {138},
      YEAR = {1999},
    NUMBER = {1},
     PAGES = {23--50},
      ISSN = {0020-9910},
   MRCLASS = {35J85 (35B65 49J40)},
  MRNUMBER = {1714335},
MRREVIEWER = {Mircea Sofonea},
       DOI = {10.1007/s002220050340},
       URL = {https://doi-org.er.lib.k-state.edu/10.1007/s002220050340},
}
\end{document}